\documentclass[journal]{IEEEtran}

\usepackage{xcolor,soul,framed, eufrak, cite} 

\colorlet{shadecolor}{yellow}
\usepackage[pdftex]{graphicx}
\usepackage{tabularx}
\usepackage{etoolbox}
\usepackage{booktabs,tabularx, colortbl}
\usepackage{multirow}
\usepackage{caption}
\usepackage{subcaption}
\usepackage{algorithm,algorithmic}
\usepackage{enumerate}
\let\labelindent\relax
\usepackage[shortlabels]{enumitem}
\def\tablename{Table}
\graphicspath{{../pdf/}{../jpeg/}}
\DeclareGraphicsExtensions{.pdf,.jpeg,.png}

\usepackage[cmex10]{amsmath}
\usepackage{amsfonts, amsmath}
\usepackage{mathtools}
\usepackage{array}
\usepackage{mdwmath}
\usepackage{mdwtab}
\usepackage{eqparbox}
\usepackage{url}
\usepackage{stmaryrd}
\usepackage{xspace}

\newcounter{propose}
\newcounter{thm}
\newcounter{define}

\newtheorem{proposition}[propose]{Proposition}
\newtheorem{definition}[define]{Definition}
\newtheorem{theorem}[thm]{Theorem}
\newtheorem{remark}{Remark}
\DeclareMathOperator*{\argmin}{arg\,min}

\newcommand{\MATLAB}{\textsc{Matlab}\xspace}
\newcommand{\cpd}[3]{\llbracket #1, #2, #3 \rrbracket}
\newcommand{\A}{{\bf A}}
\newcommand{\B}{{\bf B}}
\newcommand{\C}{{\bf C}}

\newcommand{\D}{{\bf D}}
\newcommand{\E}{{\bf E}}

\newcommand{\Ten}[1]{\mathcal{#1}}
\newcommand{\floor}[1]{\lfloor \log_2{#1}\rfloor}
\newcommand{\res}[1]{{\color{black} #1}}

\newcommand{\revision}[1]{{\color{black} #1}} 
\newcommand{\rev}[1]{{\color{black} #1}} 
\newcommand{\bc}{{\bf c}}
\newcommand{\ba}{{\bf a}}
\newcommand{\bb}{{\bf b}}

\newcommand{\bA}{{\bf A}}
\newcommand{\bB}{{\bf B}}
\newcommand{\bC}{{\bf C}}

\newcommand{\bG}{{\bf G}}

\newcommand{\bZ}{{\bf Z}}

\begin{document}
\bstctlcite{IEEEexample:BSTcontrol}
    \title{\revision{Multi-Area Distribution System State Estimation via Distributed Tensor Completion}}
     \author{
      Yajing Liu,~\IEEEmembership{Member,~IEEE,}
      Ahmed S. Zamzam,~\IEEEmembership{Member,~IEEE,}
      Andrey Bernstein,~\IEEEmembership{Member,~IEEE,}
  \thanks{The authors are with the National Renewable Energy Laboratory, Golden, CO, USA 80401 (e-mails: \texttt{\{yajing.liu, ahmed.zamzam,  andrey.bernstein\}@nrel.gov}).\newline
   This work was authored by the National Renewable Energy Laboratory, operated by Alliance for Sustainable Energy, LLC, for the U.S. Department of Energy (DOE) under Contract No. DE-AC36-08GO28308. This work was supported by the Laboratory Directed Research and Development Program at the National
Renewable Energy Laboratory.  The views expressed in the article do not necessarily represent the views of the DOE or the U.S. Government. The U.S. Government retains and the publisher, by accepting the article for publication, acknowledges that the U.S. Government retains a nonexclusive, paid-up, irrevocable, worldwide license to publish or reproduce the published form of this work, or allow others to do so, for U.S. Government purposes.}
  } 

\maketitle

\begin{abstract}
This paper proposes a model-free distribution system state estimation method based on tensor completion using canonical polyadic decomposition. In particular, we consider a setting where the network is divided into multiple areas. The measured physical quantities at buses located in the same area are processed by an area controller. A three-way tensor  is constructed to collect these measured quantities. The measurements are analyzed locally to recover the full state information of the network. A distributed closed-form iterative algorithm based on the alternating direction method of multipliers  is developed to obtain the low-rank factors of the whole network state tensor where information exchange  happens only between neighboring areas. The convergence properties of
the distributed algorithm and the sufficient conditions on the number of samples for each smaller network that guarantee the identifiability of the factors of the state tensor are presented. To demonstrate the efficacy of the proposed algorithm and to check the identifiability conditions, numerical simulations are carried out using the IEEE 123-bus system \rev{and a large-scale real utility feeder}.
\end{abstract}
\section{Introduction}
The distribution system state estimation (DSSE) task provides an estimation of voltage magnitude and angle of the network operating point (state) for distribution feeders. DSSE is a key function in the monitoring and control of electric power distribution systems \cite{Wang2019Frontiers}, such as reactive power management, demand response, distributed generation dispatch, and integration with transmission system operations \cite{McDermott2009PES}. Nowadays, as the penetration of distributed energy resources  (DERs) increases, DSSE becomes more essential to enabling normal and secure operating conditions because DERs introduce intermittent power generation and uncertain system operating conditions.

Unlike redundant measurements in transmission systems, the measured information in secondary substations of distribution systems is scarce, which poses a formidable challenge to DSSE \cite{Jabr2010IET}. To address the challenge, \revision{\cite{Dabush2021} developed a graph signal processing based weighted least squares method  and \cite{Zhang2011}  applied an optimal measurement selection approach  to the traditional Kalman filter.}
And many learning-based methods have been proposed for DSSE. For example, \cite{mestav2019bayesian} proposed a deep learning approach to Bayesian state estimation under a stochastic generation and demand model. The authors of \cite{Yang2020arXiv} built a graph neural network using several layers of unrolled Gauss-Newton iterations followed by proximal steps to account for the regularization term. The authors of \cite{zamzam2019data} presented a data-driven, learning-based neural network that can accommodate several types of measurements as well as pseudo-measurements. Also, \cite{Zamzam2020IEEEPS} proposed a novel neural network model that uses the physical structure of distribution power systems;  \cite{Ostrometzky2019arXiv} leveraged a deep neural network by incorporating physical information of the grid topology and line/shunt admittance;  \cite{donti2019matrix, Sagan2019, LiuSmartGridComm} proposed the constrained matrix completion method by combining the conventional matrix completion model with the power flow constraints; and \cite{Lin2021arXiv} developed a spatio-temporal learning approach to enhance the observability of DERs;
however, these machine learning methods require pseudo-measurements or accurate knowledge of the network model (topology of the grid or the bus admittance matrix).  Unfortunately, pseudo-measurements \revision{\cite{Filho2007}} can introduce large estimation errors \cite{Clements2011}, and accurate distribution system topology is difficult to obtain because of frequent distribution grid reconfigurations  and insufficient knowledge about the status of the network\cite{Kekatos2013IEEEPS,Muscas2015IEEIM}. As such, there is an imperative need for model-free DSSE methods.

The model-free tensor decomposition methods---including canonical polyadic decomposition (CPD), Tucker decomposition, and multilinear singular value decomposition---have been widely used in signal processing, statistics, data mining, and machine learning\rev{\cite{Shin2017IEEETKD}}, \cite{Sidiropoulos2017IEEESP, Kolda2009}. In particular, CPD has found applications in web link analysis \cite{kolda2005higher}, forecasting of the power load \cite{Song2017}, recommender systems \cite{Karatzoglou2010ACM}, \rev{PMU data recovery \cite{Osipov2020IEEETPS}}, and DSSE \cite{Zamzam2020}. In \cite{Zamzam2020}, a third-order tensor was built using the state information of the distribution network, with the form \textsc{Phase} $\times$ \textsc{Measurement} $\times$ \textsc{Time}.
\revision{
The consistency of the model generating the data is the main reason that the constructed tensor is a low-rank tensor. That is, the same power flow model governs the operation of the feeder at all time instants even if this model is unknown.
In addition, the spatial and temporal correlations in the time-series data increase the interdependency of state tensor rendering it representable using few factors in a latent space, i.e., low-rank. Thus, a model-free state estimation task is rendered a viable task using tensor completion.} The asynchronous sensor measurements across the system and the increasing complexity of the distribution system, however, call for distributed implementation of state estimation algorithms, and some distributed state estimation algorithms based on least-squares model and matrix completion were developed recently in \cite{Kekatos2013IEEEPS,zhou2019gradient,Muscas2015IEEIM,Zhu2014IEEESTSP, Sagan2019}. 

In this work, we consider a multi-area, model-free DSSE method based on tensor completion using CPD in \revision{\cite{Zamzam2020}}.  \revision{Note that the state of the network in \cite{Zamzam2020} and this work includes
the voltage phasors and  the
 active and reactive power injections at all buses/phases, which is different from the traditional state--the voltage phasors only. The way we defined the state is because when an accurate model of the network is available, the voltage phasors at all buses/phases in the network are enough to
recover any other quantity, e.g., power consumption/injection
and current flow; however, such an accurate model of
the feeder and the line parameters is often not at the system
operator’s disposal because of the change in line
parameters resulting from aging and/or unrecorded topology
changes.
} 
 
 \revision{The main contributions of this work are: \rev{
 (a) formulating the model-free multi-area distribution system state estimation as a distributed tensor completion problem, where measurements are processed locally at each area, (b) devising general identifiability conditions for the distributed tensor completion task, which represent  the minimum measuring requirements for recovering the state of the whole network, and (c) devising a distributed optimization algorithm based on ADMM that has closed-form updates for tensor low-rank factors and has a guaranteed convergent behavior. The efficacy of the proposed formulation and algorithm are demonstrated on the IEEE 123-bus system and a large-scale real utility feeder,  and the identifiability conditions are verified. 
 }
  
 }

This paper is organized as follows.
Tensor preliminaries and CPD definitions are outlined in Section II. We present the problem statement and distributed CPD model in Section III. The identifiability conditions  of  tensor completion using the multi-area CPD model  are established in Section~\ref{sec:identifiabilityconditions}.
A closed-form ADMM algorithm  is developed and its convergence is investigated in Section~\ref{sec:algorithms}. Section~\ref{sec:numericalresults} presents the simulation results, and Section~\ref{sec:conclusion} concludes the work.
\section{Tensor Preliminaries}
\label{sec:tensorpreliminaries}
In this section, we introduce some definitions related to tensors that will be later used in this manuscript. For more details, please refer to \cite{Kolda2009, Sidiropoulos2017}.

A tensor is represented as a multidimensional array. 
Consider a three-way tensor, $\Ten{X}\in\mathbb{R}^{I\times J\times K}$, with elements $\Ten{X}(i,j,k)$. The Frobenius norm of $\Ten{X}\in\mathbb{R}^{I\times J\times K}$ is defined as the square root of the sum of the
squares of all its elements, i.e., $\|\Ten{X}\|_F=\sqrt{\sum_{i=1}^I\sum_{j=1}^J\sum_{k=1}^K\Ten{X}^2(i,j,k)}$. A \emph{slab}  is a two-dimensional section of a tensor, defined by fixing one index; hence, there are three types of slabs in a three-way tensor: horizontal slabs,  $\Ten{X}(i,:,:)$,  vertical slabs, $\Ten{X}(:,j,:)$, and frontal slabs, $\Ten{X}(:,:,k)$. A fiber is a one-dimensional section of a tensor, defined by fixing the first two indices, denoted by $\Ten{X}(i,j,:)$.

A \emph{rank-one} tensor, $\Ten{X}\in\mathbb{R}^{I\times J\times K}$, is the outer product of three vectors, i.e.:
\[\mathcal{X} = {\bf a} \circ {\bf b}\circ {\bf c},\]
where ${\bf a}\in\mathbb{R}^I$, ${\bf b}\in\mathbb{R}^J$, ${\bf c}\in\mathbb{R}^K$, and the symbol ``$\circ$'' denotes the vector outer product. Any tensor can be represented as a summation of rank-one tensors in the form:
\begin{equation}
    \label{eq:CPD}
    \mathcal{X} = \sum\limits_{f=1}^F {\bf a}_f\circ {\bf b}_f \circ {\bf c}_f,
\end{equation}\vspace{-2pt}
where $F$ is a positive integer; and ${\bf a}_f \in\mathbb{R}^I, {\bf b}_f\in\mathbb{R}^J, {\bf c}_f\in\mathbb{R}^K$ for $f=1,\ldots, F$.
The minimum $F$ that can be used in (\ref{eq:CPD}) is defined as the \emph{tensor rank}, and the minimum rank polyadic decomposition is called the \emph{canonical polyadic decomposition} (CPD). 

We use matrices to denote  the tensor factors in  (\ref{eq:CPD}), i.e., $\A =[\ba_1,\ldots, \ba_F]\in\mathbb{R}^{I\times F}$, $\B =[\bb_1,\ldots, \bb_F]\in\mathbb{R}^{J\times F}$, and $\C =[\bc_1,\ldots, \bc_F]\in\mathbb{R}^{K\times F}$, \revision{in which} $\A, \B, \C$ are called the \emph{CPD factors} of the tensor $\Ten{X}$. Following \cite{Kruskal1977}, we use  the notation $\Ten{X} = \cpd{\A}{\B}{\C}$ to express the CPD equation \eqref{eq:CPD}. \revision{The authors in  \cite{Chiantini2012} derived that under some mild conditions, the CPD of $\Ten{X}$ in (\ref{eq:CPD}) is unique, and the conditions  were stated using our notations  in Theorem~1 by \cite{Zamzam2020}.}

The \emph{Kronecker product} of two matrices $\bA_1\in\mathbb{R}^{m\times n}$ and $\bB_1\in\mathbb{R}^{p\times q}$ is a matrix of size $mp\times nq$ defined by 
\begin{equation*}
       \bA_1 \otimes \bB_1 =\begin{bmatrix}
  a_{11}\bB_1 &\cdots &a_{_{1n}}\bB_1\\
    \vdots &\ddots &\vdots\\
    a_{_{m1}}\bB_1 &\cdots &a_{_{mn}}\bB_1
    \end{bmatrix}.
\end{equation*}
The \emph{Khatri-Rao product} is the ``matching columnwise" Kronecker product. The Khatri-Rao product of 
$\bA\in\mathbb{R}^{I\times F}$ and $\bB\in\mathbb{R}^{J\times F}$ is a matrix of size $IJ \times F$ and defined by
\begin{equation*}
       \bA \odot \bB =\begin{bmatrix}
  \ba_1\otimes \bb_1 &\ba_2\otimes \bb_2& \cdots &\ba_F\otimes \bb_F
    \end{bmatrix}.
\end{equation*}

\section{Problem Formulation}
\label{sec:distributedprobform}
In \revision{a general multi-phase} electric power distribution system, we define the state of the network as the collection of active and reactive power injections as well as the voltage phasors at \revision{all the phases of all the buses. In the following, we refer to the pair (bus, phase) as simply ``phase'' for brevity.} The state estimation task aims to estimate the unknown values of the state variables using available measurements.
In this work, we construct a three-way tensor, $\Ten{X}\in\mathbb{R}^{I\times J\times K}$, that collects nodal physical quantities (states) for multiple time instants. The tensor $\Ten{X}$ is in the form \textsc{Phase} $\times$ \textsc{Measurement} $\times$ \textsc{Time}. That is, the element $\Ten{X}(i,j,k)$ represents the measurement $j$ for phase $i$ taken at time $k$, \revision{and $I$ is the number of phases for all the buses in the network, $J$ is the number of the measurement types, and $K$ is the number of time steps}. \revision{The measurement types considered in this paper include the real part of the voltage phasor, imaginary part of the voltage phasor, voltage magnitude, active power, and reactive power, so $J=5$.  We  remark  that  our  model  is  able  to  accommodate any measurement known to have some correlation with other measurements in the tensor.}
 The low-rank CPD representation of the tensor $\Ten{X}$ can be expressed as follows:
\begin{align}
\Ten{X}(i, j, k) = \sum_{f=1}^{F} \A(i, f) \B(j, f) \C(k, f),
\end{align}
where $F$ denotes the rank of the tensor, and the matrices $\A$, $\B$, and $\C$ denote the low-rank factors.  

In general, the state tensor $\Ten{X}$ has a low-rank structure for three reasons: 1) the  temporal correlation that exists as a result of similar conditions affecting power consumption at all phases; 2) the spatial correlation between different locations in a power grid affecting solar consumption, for example;  and  3) the  relationship  between  different  types  of  measurements via the power flow model, which can be well-approximated using linear flow models. Therefore, this motivates us to apply the following CPD model \cite{Zamzam2020} to recover the unknown values of $\Ten{X}$:
 \begin{align}\label{eq:centralizedmodel}
    \min_{\A, \B, \C} \frac{1}{2}\| \Ten{W} * \big(\Ten{X} - {\cpd{\A}{\B}{\C}} \big)\|_F^2,
\end{align}
where 
$\Ten{X}\in\mathbb{R}^{I\times J\times K}$ is the tensor to be completed; ${\A},{\B},{\C}$ are the CPD factors of $\Ten{X}$ with
$\bA\in\mathbb{R}^{I\times F}, \bB\in\mathbb{R}^{J\times F},\bC\in\mathbb{R}^{K\times F}$; 
$\Ten{W}\in \{0,1\}^{I\times J\times K}$ is a binary observation tensor with 1 for every known element, 0 otherwise; \revision{and $*$ is the Hadamard or element-wise product. Note that we use the same notation with \cite{tensorlab} in (\ref{eq:centralizedmodel}).}

In \cite{Zamzam2020}, the centralized model-free state estimation methods for solving (\ref{eq:centralizedmodel}) have been proposed. Unfortunately, these methods require a high degree of synchronization among measuring units across the system, which is challenging to achieve, especially in large-scale distribution feeders; thus, we aim to build a multi-area CPD computation model and develop a distributed state estimation method for solving the proposed model, which circumvents the issue of requiring network-wide measurement synchronization. In addition, it will  also alleviate computational burdens for the control center. The distributed state estimation method does not need a centralized coordination system operator but only a communication network for exchanging the boundary and consensus variables among adjacent regions. The multi-area CPD model of (\ref{eq:centralizedmodel}) is built as follows.

Consider a general three-phase distribution network, where $\Ten{X}\in\mathbb{R}^{I\times J\times K}$ denotes the state tensor. Assume that the distribution system is partitioned into $N$ areas along the first dimension (phases). Let $\mathcal{I}_i$ denote the set of phases in the $i$-th area, where $I_i$ denotes the number of phases in area $i$, and let $\mathcal{N}(i)$ denote the areas adjacent to area $i$. Note that this adjacency is defined over the communications network between area controllers that are not necessarily conforming with the electric network because we assume that the electric network topology is not known. We also assume that the set $\mathcal{I}_i$ is not empty for all $i$.
Let $\Ten{X}_i$, $\Ten{W}_i$, and $\A_i$ be the corresponding parts of $\Ten{X}$, $\Ten{W}$, and $\A$ with respect to area $i$, which satisfy that $\Ten{X} = [\Ten{X}_1; \ldots; \Ten{X}_N]$, $\Ten{W} = [\Ten{W}_1; \ldots; \Ten{W}_N]$, and $\A=[\A_1;\ldots; \A_N]$. Each area $i$ is required to communicate with its adjacent areas to exchange the consensus variables $\B_i$ and $\C_i$; thus,
for each area $i$,  $\B_i$ and $\C_i$ are defined as the basis matrices satisfying $\B_i = \B$ and $\C_i = \C$. \revision{Note that, the communication between neighboring areas in the feeder does not require high-level synchronization and several asynchronous optimization approaches can be adopted.} Then (\ref{eq:centralizedmodel}) can be written in a distributed form as follows:
\begin{equation}
   \label{eq:distributedCPDmodel}
    \begin{aligned}
   \min_{\{\A_i\}, \{\B_i\}, \{\C_i\}} &\sum_{i=1}^N\frac{1}{2}\| \Ten{W}_i * \big(\Ten{X}_i - {\cpd{\A_i}{\B_i}{\C_i}} \big)\|_F^2,\\
     \ \text{s.t.}\ &\B_i = \B_{j}, j\in\mathcal{N}(i),\\
     &  \C_i = \C_j, j\in\mathcal{N}(i).
    \end{aligned}
\end{equation}

Auxiliary variables are introduced to enforce the connected areas to consent on the basis matrices. We introduce $\B_{ij}= \B_i = \B_j$ and 
$\C_{ij}= \C_i = \C_j$ if $j\in\mathcal{N}(i)$ to guarantee that the equality constraints in (\ref{eq:distributedCPDmodel}) hold. Including the new variables $\{\B_{ij}\}$ and $\{\C_{ij}\}$, (\ref{eq:distributedCPDmodel}) becomes:
\begin{equation}
   \label{eq:auxidistributedCPDmodel}
    \begin{aligned}
    \min_{\substack{\{\A_i\}, \{\B_i\}, \{\C_i\}\\ \{\B_{ij}\}, \{\C_{ij}\}}} &\sum_{i=1}^N\frac{1}{2}\| \Ten{W}_i * \big(\Ten{X}_i - {\cpd{\A_i}{\B_i}{\C_i}} \big)\|_F^2,\\
     \ \text{s.t.}\ &\B_i = \B_{ij}=\B_{j}, j\in\mathcal{N}(i),\\
     &  \C_i = \C_{ij}=\C_{j}, j\in\mathcal{N}(i).
    \end{aligned}
\end{equation}
The following section establishes the identifiability conditions of the tensor completion using the multi-area CPD model~(\ref{eq:auxidistributedCPDmodel}).

\section{Identifiability Conditions}
\label{sec:identifiabilityconditions}
The identifiability conditions of the  tensor
completion using the centralized CPD model were studied by  \cite{Zamzam2020} and \cite{Kanatsoulis2020}; however, the  identifiability conditions in these two references require  all the measurements to be synchronized across the entire system, and the conditions 
 apply only to specific uniform sampling schemes. 
In this section, we provide more general identifiability conditions  of the tensor completion using the multi-area CPD model (\ref{eq:distributedCPDmodel}), which also covers the case of the centralized CPD model. Our identifiability conditions  allow for asynchronous samplings among different areas, and they apply to any sampling scheme.
We first introduce some sampling notations and schemes, and then we provide the general identifiability conditions.
\subsection{Sampling Notations and Schemes}
\label{subsect:samplingnotations}
Assume that the distribution system is partitioned into $N$ areas. Review the notations in Section~\ref{sec:tensorpreliminaries} and Section~\ref{sec:distributedprobform}, $\mathcal{X}={\cpd{\A}{\B}{\C}}\in\mathbb{R}^{I\times J\times K}$ denotes the whole tensor to be completed, and $\mathcal{X}_n(\mathcal{I}_n, :,:)={\cpd{\A_n}{\B_n}{\C_n}}\in\mathbb{R}^{I_n\times J\times K}$ denotes the subtensor to be completed in area $n$, satisfying that $\Ten{X} = [\Ten{X}_1; \ldots; \Ten{X}_N]$, $\A = [\A_1; \ldots; A_N]$, $\B=\B_n$, and $\C=\C_n$, where $\mathcal{I}_n$ is the set of phases in area $n$. Without loss of generality, let $\mathcal{I}_n=\{{\sum_{m=1}^{n-1}I_m+1,\ldots, \sum_{m=1}^{n}I_m}\}$ with $\sum_{n=1}^N I_n=I$.

In area $n$, let $\Ten{\bar{X}}_n$ denote the sampled subtensor. Let  $\mathcal{P}_n\subseteq \mathcal{I}_n$, $\mathcal{M}_n\subseteq\{1,\ldots, J\}$, and $\mathcal{T}_n\subseteq\{1,\ldots, K\}$ be the sets of identifiable rows of $\A_n$, $\B_n$, and $\C_n$ through the available samples of $\Ten{\bar{X}}_i$, respectively. Using  the sampled subtensor $\Ten{\bar{X}}_n$ in area $n$, we recover $\A_n$ and partial rows of $\B$ and $\C$. Our aim is to recover the whole tensor $\Ten{X}$ through the sampled subtensors in all areas. In the following proposition, we  provide the conditions under which the whole tensor $\mathcal{X}$ can be recovered from the sampled subtensors.

\subsection{Main Result: Identifiability}
Before we state the main result, we introduce some definitions that will prove useful later on.

\begin{definition}
Area $n$ and area $m$ are mutually identifiable if i) $|\mathcal{M}_n\cap \mathcal{M}_m|\geq 2 \ \text{and}\  |\mathcal{T}_n\cap \mathcal{T}_m|\geq 1$ hold, or \\
ii) $|\mathcal{T}_n\cap \mathcal{T}_m|\geq 2 \ \text{and}\  |\mathcal{M}_n\cap \mathcal{M}_m|\geq 1$ hold.
\end{definition}

\revision{The concept of mutually identifiable partitions is introduced in this work to facilitate the presentation of the main theoretical result. Since the sampled measurements and time instants differ between areas of the feeder, we denote the set of measurements and time steps whose low-rank representation is identifiable from the sampled measurements at area $n$ by $\mathcal{M}_n$ and $\mathcal{T}_n$, respectively. Then, two areas are \emph{mutually identifiable} if their identifiable components have an intersection as detailed in the result. This intersection guarantees that identified factors from the two mutually identifiable partitions can be combined, i.e., their respective rotation and scaling ambiguities can be resolved.}

\begin{definition}
Consider each area of the distribution system as a vertex; if two vertices are connected, we say that there is an edge between the two vertices.
Define a graph $\mathcal{G}(\mathcal{V}, \mathcal{E})$, where $\mathcal{V}$ denotes the set of $N$ vertices, and $\mathcal{E}$ denotes the set of all edges connecting any two vertices. 
Then, we say that vertices representing area $n$ and area $m$ are connected in the graph $\mathcal{G}$ if and only if the two areas are mutually identifiable.
\end{definition}

\revision{Note that the graph considered here is independent of the physical systems graph. Since the proposed approach is model-free, the measurement requirements used are not function of the electric feeder topology. Additionally, the quality of the estimates provided from the proposed approach should not be dependant on the partitioning of the system as long as the conditions presented next are satisfied. Nonetheless, utilizing the feeder topology to design partitioning strategies may lead to enhancement in convergence speed of ADMM-based decentralized optimization algorithm.}

\begin{proposition}
\label{proposition:identifiability}
Assume that  $\A_n, \B_n$, and $\C_n$ are drawn
from some joint absolutely continuous distribution with respect
to the Lebesgue measure in $\mathbb{R}^{(I_n+J+K)F}$. If the sampled subtensors $\Ten{\bar{X}}_n$, for every area $n$ ($n=1,\ldots, N$), admit essentially unique CPDs 
and the following conditions are satisfied:
	\begin{enumerate}[start=1,label={ \arabic*)}]
		\item $ \mathcal{P}_n=\mathcal{I}_n\  \forall n=1,\ldots, N$, $\bigcup_{n=1}^N \mathcal{M}_n=\{1,\ldots,J\}$, and  $\bigcup_{n=1}^N \mathcal{T}_n=\{1,\ldots, K\}$,
		\item The graph $\mathcal{G}(\mathcal{V}, \mathcal{E})$ is connected,
	\end{enumerate}
	then $\Ten{X}$ can be recovered almost surely from the observed elements.
\end{proposition}
The proof of this proposition is relegated to Appendix A.
\begin{figure}
    \centering
  \includegraphics[scale=0.5,trim= 190 130 20 20, clip]{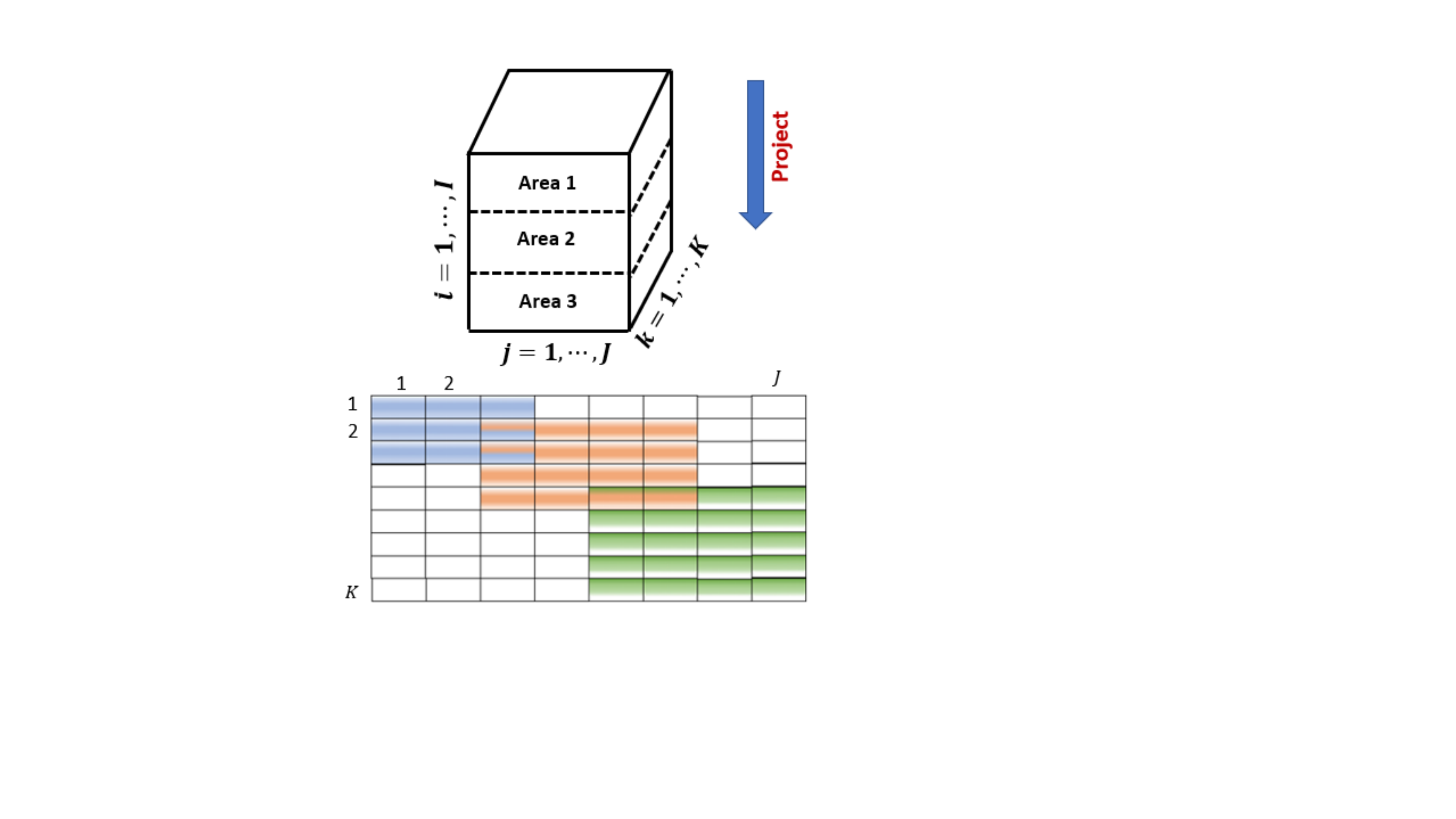}
    \caption{Illustration of connected partitioning}
    \label{fig:illustration}\vspace{-15pt}
\end{figure}
\begin{remark}
Fig.~\ref{fig:illustration} illustrates the definition of the graph $\mathcal{G}(\mathcal{V}, \mathcal{E})$ being connected. The system is partitioned into three areas based on the first dimension. The blue area, orange area, and green area denote the sets of identifiable rows of $\B_n$ and $\C_n$ ($n=1,2,3$). The figure shows that $|\mathcal{M}_1\cap \mathcal{M}_2|=1, |\mathcal{T}_1\cap \mathcal{T}_2|=2$, and $|\mathcal{M}_2\cap \mathcal{M}_3|=2, |\mathcal{T}_2\cap \mathcal{T}_3|=1$, which implies that Area 1 and Area 2 are mutually identifiable, and Area 2 and Area 3 are mutually identifiable; thus, the graph with the three areas being vertices is connected. 
\end{remark}

\section{Algorithms}
\label{sec:algorithms}
\subsection{ADMM Algorithm}
\label{section:disADMM}
In this section, we develop a distributed algorithm based on ADMM that enjoys closed-form updates to solve the multi-area CPD model~(\ref{eq:auxidistributedCPDmodel}). \revision{The updates in this approach are done locally at each area after receiving updates from neighboring areas. On the other hand, the updated variables are sent to neighboring areas following the local updates. This allows for asynchronous distributed optimization to be deployed. This relaxes the expensive requirement of having high-level synchronization between all areas that is often required in centralized optimization methods.}

Let ${\bf{\Gamma}}_{ij}$ and ${\bf{\Lambda}}_{ij}$ denote the Lagrange multipliers associated with the pair of equation constraints in (\ref{eq:auxidistributedCPDmodel}), respectively.
For simplicity, we use ${\bf{\Gamma}}_{ij}$ and ${\bf{\Lambda}}_{ij}$ to denote the scaled forms of the dual variables, i.e., ${\bf{\Gamma}}_{ij}/\mu$ and ${\bf{\Lambda}}_{ij}/\lambda$.  
The scaled augmented Lagrangian function (please refer to \cite{Boyd_admm} for details) of (\ref{eq:auxidistributedCPDmodel}) can be written as follows:
\begin{equation}
     \label{eq:augmentedlagrangianfunc}
     \begin{aligned}
    \mathcal{L}(&\{\A_i\}, \{\B_i\}, \{\C_i\}, \{\B_{ij}\}, \{\C_{ij}\}, \{{\bf{\Gamma}}_{ij}\}, \{{\bf{\Lambda}}_{ij}\})\\
     \ = &\sum_{i=1}^N\frac{1}{2}\| \Ten{W}_i * \big(\Ten{X}_i - {\cpd{\A_i}{\B_i}{\C_i}} \big)\|_F^2\\
    &+\sum_{i=1}^N\sum_{j\in\mathcal{N}(i)}\frac{\mu}{2}\|\B_{i}-\B_{ij}+{\bf{\Gamma}}_{ij}\|_F^2-\frac{\mu}{2}\|{\bf{\Gamma}}_{ij}\|_F^2\\
    &+\sum_{i=1}^N\sum_{j\in\mathcal{N}(i)}\frac{\lambda}{2}\|\C_{i}-\C_{ij}+{\bf{\Lambda}}_{ij}\|_F^2-\frac{\lambda}{2}\|{\bf{\Lambda}}_{ij}\|_F^2.
    \end{aligned}
 \end{equation}

Using the scaled dual variables, the ADMM updates are given \revision{as follows.}
\newline The primal variable updates are: 
\[\A_i^{k+1} = \argmin_{\A_i} \frac{1}{2}\| \Ten{W}_i * \big(\Ten{X}_i - {\cpd{\A_i}{\B_i^k}{\C_i^k}} \big)\|_F^2,
\]
\begin{align*}
    \B_i^{k+1} = &\argmin_{\B_i} \left\{\frac{1}{2}\| \Ten{W}_i * \big(\Ten{X}_i - {\cpd{\A_i^{k+1}}{\B_i}{\C_i^k}} \big)\|_F^2\right.\\
    &\ \left.+\sum_{j\in\mathcal{N}(i)}\frac{\mu}{2}||\B_i-\B_{ij}^k+{\bf{\Gamma}}_{ij}^{k}||_F^2\right\},
\end{align*}
\begin{align*}
    \C_i^{k+1} = &\argmin_{\C_i} \left\{\frac{1}{2}\| \Ten{W}_i * \big(\Ten{X}_i - {\cpd{\A_i^{k+1}}{\B_i^{k+1}}{\C_i}} \big)\|_F^2\right.\\
    &\ \left.+\sum_{j\in\mathcal{N}(i)}\frac{\lambda}{2}||\C_i-\C_{ij}^k+{\bf{\Lambda}}_{ij}^{k}||_F^2\right\}.
\end{align*}
The auxiliary variables updates are: 
\[\B_{ij}^{k+1} = \frac{\B_i^{k+1}+\B_j^{k+1}}{2}+\frac{{\bf{\Gamma}}_{ij}^k+{\bf{\Gamma}}_{ji}^k}{2},\]
\[\C_{ij}^{k+1} = \frac{\C_i^{k+1}+\C_j^{k+1}}{2}+\frac{{\bf{\Lambda}}_{ij}^k+{\bf{\Lambda}}_{ji}^k}{2}.\]
The dual variable updates are:
\begin{align*}
    {\bf{\Gamma}}_{ij}^{k+1} = {\bf{\Gamma}}_{ij}^k +(\B_i^{k+1}-\B_{ij}^{k+1}),
\end{align*}
\begin{align*}
    {\bf{\Lambda}}_{ij}^{k+1} = {\bf{\Lambda}}_{ij}^k +(\C_i^{k+1}-\C_{ij}^{k+1}).
\end{align*}
From the update rule of ${\bf{\Gamma}}_{ij}$, we have: 
\[{\bf{\Gamma}}_{ij}^{k+1} + {\bf{\Gamma}}_{ji}^{k+1}= {\bf{\Gamma}}_{ij}^k+{\bf{\Gamma}}_{ji}^{k} +(\B_i^{k+1}+\B_j^{k+1})-2\B_{ij}^{k+1}.\]
Substituting the update rule of $\B_{ij}$ into the above equation results in 
${\bf{\Gamma}}_{ij}^{k+1} + {\bf{\Gamma}}_{ji}^{k+1}=0,$
which implies that: \[\B_{ij}^{k} = \frac{\B_i^{k}+\B_j^{k}}{2}.\]
Similarly, we have: \[\C_{ij}^{k} = \frac{\C_i^{k}+\C_j^{k}}{2}.\]

Substituting the above two equations into the primal and dual updates, we have the following updates for the primal and dual variables:
\begin{align*}
    \A_i^{k+1} &= \argmin_{\A_i}\frac{1}{2} \| \Ten{W}_i * \big(\Ten{X}_i - {\cpd{\A_i}{\B_i^k}{\C_i^k}} \big)\|_F^2,
\end{align*}
\begin{align*}
    \B_i^{k+1}& = \argmin_{\B_i} \left\{\frac{1}{2}\| \Ten{W}_i * \big(\Ten{X}_i - {\cpd{\A_i^{k+1}}{\B_i}{\C_i^k}} \big)\|_F^2\right.\\
    &\ \ \  \left.+\sum_{j\in\mathcal{N}(i)}\frac{\mu}{2}\|\B_i-{(\B_i^{k}+\B_j^{k})}/{2}+{\bf{\Gamma}}_{ij}^{k}\|_F^2\right\},
\end{align*}
\begin{align*}
    \C_i^{k+1} &= \argmin_{\C_i} \left\{\frac{1}{2}\| \Ten{W}_i * \big(\Ten{X}_i - {\cpd{\A_i^{k+1}}{\B_i^{k+1}}{\C_i}} \big)\|_F^2\right.\\
    &\ \ \  \left.+\sum_{j\in\mathcal{N}(i)}\frac{\lambda}{2}{\|}\C_i-(\C_i^{k}+\C_j^{k})/{2}+{\bf{\Lambda}}_{ij}^{k}{\|}_F^2\right\},
\end{align*}
\begin{equation}
\label{eq:dual1update}
    {\bf{\Gamma}}_{ij}^{k+1} = {\bf{\Gamma}}_{ij}^k +{(\B_i^{k+1}-\B_j^{k+1})}/{2},
j\in\mathcal{N}(i)
\end{equation}
\begin{equation}
    \label{eq:dual2update}
    {\bf{\Lambda}}_{ij}^{k+1} = {\bf{\Lambda}}_{ij}^k +{(\C_i^{k+1}-\C_j^{k+1})}/{2}, j\in\mathcal{N}(i).
\end{equation}
Let $\Ten{W}_i^{(j)}$, $\Ten{X}_i^{(j)}$ for $j=1,2,3$ denote the mode-$j$ unfoldings of tensors $\Ten{W}_i$ and $\Ten{X}_i$, respectively. Define $\underline{\bf w}_i^{(j)}=\text{vec}(\Ten{W}_i^{(j)})$, $\underline{\bf x}_i^{(j)}=\text{vec}(\Ten{X}_i^{(j)})$, ${\bf a}_i = \text{vec}(\A_i)$, ${\bf b}_i = \text{vec}(\B_i)$,  and ${\bf c}_i = \text{vec}(\C_i)$. Then finding the primal variables $\A_i$, $\B_i$, and $\C_i$ is equal to finding ${\bf a}_i$, ${\bf b}_i$, and ${\bf c}_i$, respectively.
\revision{Using the matricization of  tensors and the vectorization of matrices, by the first order necessary condition,}
the closed-form updates for the primal variables 
${\bf a}_i$, ${\bf b}_i$, and ${\bf c}_i$ are derived 
as follows:
\begin{equation}
\label{eq:primal1update}
    {\bf a}_i^{k+1} 
    =(\D_i^{kT}\D_i^k)^{-1}\D_i^{kT}(\underline{\bf w}_i^{(1)}* \underline{\bf x}_i^{(1)}),
\end{equation}
where $\D_i^k=\text{diag}(\underline{\bf w}_i^{(1)})((\C_i^k\odot \B_i^k)\otimes {\bf I}_{_{I_i\times I_i}})$.
\begin{align}
\label{eq:primal2update}
    {\bf b}_i^{k+1}
    &= (\E_i^{kT}\E_i^k+\mu\sum_{j\in\mathcal{N}(i)}{\bf I}_{{FJ\times FJ}})^{-1}(\E_i^{kT} (\underline{\bf w}_i^{(2)}* \underline{\bf x}_i^{(2)})+\nonumber\\
    &\quad\mu\sum_{j\in\mathcal{N}(i)}\text{vec}({(\B_i^{k}+\B_j^{k})}/{2}-{\bf{\Gamma}}_{ij}^{k})),
\end{align}
where $\E_i^k  =\text{diag}(\underline{\bf w}_i^{(2)})((\C_i^k\odot \A_i^{k+1})\otimes {\bf I}_{_{J\times J}})$.
\begin{align}
\label{eq:primal3update}
    {\bf c}_i^{k+1}
    &=(\bG_i^{kT}\bG_i^k+\lambda\sum_{j\in\mathcal{N}(i)}{\bf I}_{{FK\times FK}})^{-1}(\bG_i^{kT}(\underline{\bf w}_i^{(3)} * \underline{\bf x}_i^{(3)})+\nonumber\\
    &\quad\lambda\sum_{j\in\mathcal{N}(i)}\text{vec}((\C_i^{k}\C_j^{k})/{2}-{\bf{\Lambda}}_{ij}^{k})),
\end{align}
where  $\bG_i^k =\text{diag}(\underline{\bf w}_i^{(3)})((\B_i^{k+1}\odot \A_i^{k+1})\otimes {\bf I}_{_{K\times K}})$. \revision{The readers can refer to \cite{Liavas2015} for deriving the closed-form updates for the primal variables in a similar model.}

The closed-form ADMM algorithm for the multi-area CPD model (\ref{eq:auxidistributedCPDmodel}) is given as follows:
\begin{algorithm}
\caption{ADMM Algorithm for Distributed CPD Model}
 \begin{algorithmic}[1]
 \renewcommand{\algorithmicensure}{\textbf{Initialize:}}
\ENSURE Obtain the CPD factors $\A_i^1, \B_i^1, \C_i^1$ of the sub-measurement tensor $\Ten{W}_i *\Ten{X}_i$ for $i=1,\cdots, N$,  then set ${\bf a}_i^1, {\bf b}_i^1, {\bf c}_i^1$ as the vector form of $\A_i^1, \B_i^1, \C_i^1$, respectively.  Let ${\bf{\Gamma}}_{ij}^1=0, {\bf{\Lambda}}_{ij}^1=0$ for $j\in\mathcal{N}(i)$.
  \FOR {$k = 1,\ldots,n$ 
  }
\STATE 
Use (\ref{eq:primal1update}), (\ref{eq:primal2update}), (\ref{eq:primal3update}) to update ${\bf a}_i^k, {\bf b}_i^k, {\bf c}_i^k$, respectively.

Update ${\bf{\Gamma}}_{ij}^k, {\bf{\Lambda}}_{ij}^k$ for $j\in\mathcal{N}(i)$ using (\ref{eq:dual1update}) and (\ref{eq:dual2update}), respectively.
  \ENDFOR
 \RETURN $\Ten{X}=$ 
 
 $[\cpd{\A_1^{n+1}}{\B_1^{n+1}}{\C_1^{n+1}}; \cdots; \cpd{\A_N^{n+1}}{\B_N^{n+1}}{\C_N^{n+1}}]$ 
 \end{algorithmic} 
 \label{Agorthim:ADMM}
\end{algorithm}

\subsection{Convergence}
In this section, we discuss the convergence of the distributed ADMM algorithm devised in the previous section.

We say a point $\bZ=\left(\{{\bf a}_i\}, \{{\bf b}_i\}, \{{\bf c}_i\},\B_{ij},\C_{ij}, \{{\bf{\Gamma}}_{ij}\}, \{{\bf{\Lambda}}_{ij}\}\right)$ satisfies the Karush-Kuhn-Tucker (KKT) conditions for problem (\ref{eq:auxidistributedCPDmodel}) if there exists $\bZ$ such that: 
 \begin{subequations}
 \label{eq:KKTconditions}
      \begin{equation}
     \label{eq:stationarityA}
       \D_i^T(\underline{\bf w}_i^{(1)}* \underline{\bf x}_i^{(1)}-\D_i{\bf a}_i)=0,
 \end{equation}
  \begin{equation}
     \label{eq:stationarityB}
       \begin{aligned}
    &(\E_i^T\E_i+\mu\sum_{j\in\mathcal{N}(i)}{\bf I}_{{FJ\times FJ}}) {\bf b}_i=\\
    &\quad\quad\quad\quad\E_i^T (\underline{\bf w}_i^{(2)}* \underline{\bf x}_i^{(2)})+\mu\sum_{j\in\mathcal{N}(i)}\text{vec}(\B_{ij}-{\bf{\Gamma}}_{ij}),
\end{aligned}
 \end{equation}
  \begin{equation}
     \label{eq:stationarityC}
       \begin{aligned}
    &(\bG_i^T\bG_i+\lambda\sum_{j\in\mathcal{N}(i)}{\bf I}_{{FK\times FK}}) {\bf c}_i=\\
    &\quad\quad\quad\bG_i^T (\underline{\bf w}_i^{(3)}* \underline{\bf x}_i^{(3)})+\lambda\sum_{j\in\mathcal{N}(i)}\text{vec}(\C_{ij}-{\bf{\Lambda}}_{ij}),
\end{aligned}
 \end{equation}
 \begin{equation}
     \label{eq:primalfeasibility1}
     \B_i -\B_{j}=0, j\in\mathcal{N}(i),
 \end{equation}
  \begin{equation}
     \label{eq:primalfeasibility2}
     \C_i -\C_{j}=0, j\in\mathcal{N}(i).
 \end{equation}
 \end{subequations}
 
 \begin{theorem}
 If the multiplier sequence $(\{{\bf{\Gamma}}_{ij}^k\}, \{{\bf{\Lambda}}_{ij}^k\})$ for $i=1,\ldots, N, j\in \mathcal{N}(i)$ is bounded and satisfies:
\begin{equation}
\label{ineq:boundassumption}
    \sum\limits_{k=0}^{\infty}\sum_{i=1}^N\sum_{j\in\mathcal{N}(i)}\|{\bf{\Gamma}}_{ij}^{k+1}-{\bf{\Gamma}}_{ij}^k\|_F^2+\|{\bf{\Lambda}}_{ij}^{k+1}-{\bf{\Lambda}}_{ij}^k\|_F^2<\infty,
\end{equation}
then any accumulation point of $\{\bZ^k\}$ generated by the ADMM for problem 
(\ref{eq:auxidistributedCPDmodel}) is a KKT point satisfying (\ref{eq:KKTconditions}).
 \end{theorem}
\begin{proof}
See Appendix A.
\end{proof}
\section{Numerical Results}
\label{sec:numericalresults}
In this section, we demonstrate the performance of the proposed multi-area model-free distribution system state estimation approach on the IEEE 123-bus system \rev{and a real utility feeder}, and we compare that with the centralized SDF-NLS solver used in \cite{Zamzam2020}. \rev{The data
we use for our implementation were simulated at 1-minute
resolution using power flow analysis with diversified load and
solar profiles created for each bus using real solar irradiance
and real load consumption data. For both the IEEE 123-bus system and real utility feeder, we assumed
a constant power factor loads, but when renewable energy injections are subtracted (added with
a different sign to obtain the net power injections), the net power injections  do not have a constant
power factor. Thus, the active and reactive power injections are not completely dependant if only
the net injections are observed.}

\revision{\rev{ The measurement types considered in the simulation include the real/imaginary parts of the voltage phasor, voltage magnitude, and the net active/reactive power.} \rev{The first two measurement types are assumed to be obtained from the very limited PMU measuring devices installed within the the distribution feeder, and the last three measurements are assumed to be obtained from smart meters of SCADA systems.}  The tensor is constructed using 72 time steps that are sampled  every 20 minutes, thus this nonconsecutive data contain the state information
of the network over 24 hours and the tensor is of size $263\times 5\times 72$.} \rev{Note that, we have experimented using the angles at each phase as a separate measurement. However, our simulations indicated that the rank required to approximate the tensor increases significantly when this measurement is added. This is mainly due to the angles being in three different ranges, i.e., $0$, $\frac{2\pi}{3}$, and $\frac{-2\pi}{3}$. Thus, we have excluded the angles from the state variables and considered it as an implicit state variable from the first stages of this work.}

The known measurements are added with noise with zero mean and 1\% of the true value as the  standard deviation.  The  mean absolute percentage error (MAPE) for the voltage magnitude ($|V|$) and mean absolute error (MAE) for the voltage angle ($\theta$), active power ($P$), and reactive power ($Q$) are used to measure the performance of our closed-form decentralized approach and the centralized SDF-NLS solver. \revision{And the MAPE and MAE are defined as:
$\text{MAPE}: = \frac{1}{n}\sum\frac{|\text{Actual-Estimated}|}{\text{Actual}}\times 100\%$ and $\text{MAE}: = \frac{1}{n}\sum(|\text{Actual-Estimated}|)$, where $n$ is the total number of elements. } The sampling schemes considered are slab sampling, fiber sampling, and the mixed sampling, and the definitions of slabs and fibers were introduced in Section~\ref{sec:tensorpreliminaries}. \revision{The experiments below are performed on a laptop with 1.9 GHZ CPU and 32 GB RAM. \rev{For our simulation experiments, at each phase, we normalize each type of its known data by subtracting its mean over
all the time steps, and recover back finally, so the errors are calculated based on its original data.}}

\begin{figure}
    \centering
  \includegraphics[scale=0.5,trim= 80 220 20 250, clip]{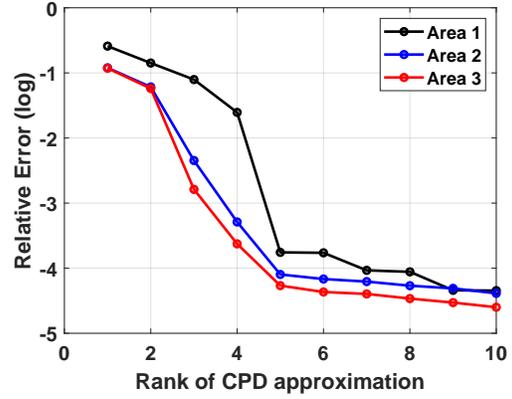}
    \caption{Relative error of the computed CPD approximation}
    \label{fig:lowrank}\vspace{-15pt}
\end{figure}

The system is partitioned into three areas with the areas comprising 87 phases, 87 phases, and 89 phases, respectively. Therefore, the three subtensors $\mathcal{X}_1$, $\mathcal{X}_2$, and $\mathcal{X}_3$ are of sizes $87\times 5\times 72$, $87\times 5\times 72$, and $89\times 5\times 72$, respectively. By the identifiablity condition
in \cite{Chiantini2012}, if the rank ($F$) of the subtensors satisfies $F\leq 2^{\floor{5}+\floor{72}-2}=32$, then the low-rank factors are essentially unique almost surely.
Fig.~\ref{fig:lowrank} shows that when the rank of the CPD approximation for each subtensor is equal to seven, the relative error between the subtensor  and  its  rank-$k$ CPD approximation is  less than $10^{-4}$, which implies that the three  subtensors can be well approximated using low-rank tensor decompositions. The relative error between  the  subtensor  and  its  rank-$k$ CPD approximation is  defined as the ratio of the Frobenius norm square of the residual between the subtensor and its rank-$k$ CPD approximation and the Frobenius norm square of the subtensor.

\rev{The authors \cite{ZhangIEEETSP2018} proposed a model-free multichannel Hankel matrix completion method for missing data recovery and provided performance guarantees. However, the matrix can not be recovered if one entire row (which comprises one measurement type for all locations at different time steps) is missing.}
\rev{The proposed approach is a distributed fashion of the centralized model \cite{Zamzam2020}, and allows a entire row missing of one slice of the tensor to be recovered. Thus, the proposed method will be compared with the centralized approach \cite{Zamzam2020} and a model-based decentralized matrix completion method \cite{Sagan2019}. The simulation results showed the efficacy of the distributed approach and that the model-free approach can achieve close estimation accuracy with model-based methods.}
\subsection{Slab Sampling}
\label{sec:slabsampling}
\begin{figure}
    \centering
  \includegraphics[scale=0.5]{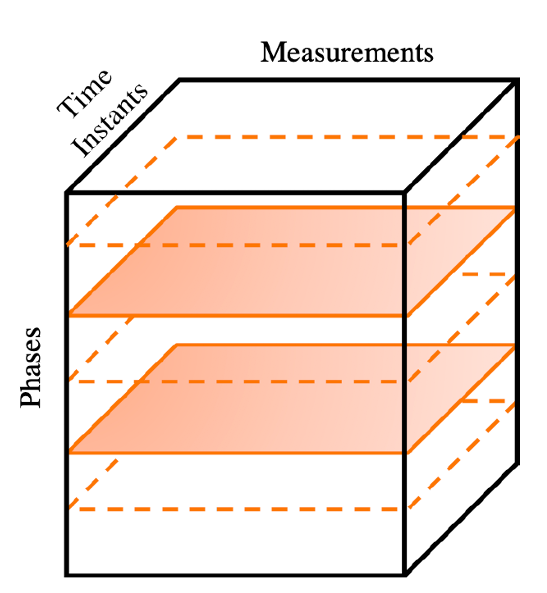}
    \caption{Horizontal slab sampling example.}
    \label{fig:slab}
\end{figure}

For slab sampling, we consider horizontal and frontal slab sampling schemes for each area. We use $\Ten{\bar{X}}_n(\mathcal{S}_n^p, :,:)$  to denote the horizontal slab sampled in area $n$, where the phases with all measurements available at all time instances are collected in the set $\mathcal{S}_n^p\subseteq\mathcal{I}_n$; Fig. \ref{fig:slab} shows an example of horizontal slab sampling. We use $\Ten{\bar{X}}_n(:,:, \mathcal{S}_n^t)$ to denote the frontal slab sampled in area $n$, where all measurements for all phases are sampled at the time steps collected in $\mathcal{S}_n^t\subseteq\{1,\ldots, K\}$.  \rev{It is worth noting that, the slab sampling measuring scheme requires frontal slabs to be sampled from the state tensor which includes the real and imaginary parts of voltage phasors at all phases. This is impractical due to the limited availability of such devices in distribution networks. However, the proposed distributed framework allows for requiring such frontal slab to be sampled at exactly one partition of the network according to Theorem 1. Thus, the operator may assume that one partition containing all phases with PMU measurements has complete frontal slab sampling while the other network partitions do not require this measurement type.}

 For horizontal slab samplings, we sample eight phases randomly from nonzero load phases for each area, so $|\mathcal{S}_n^p|=8$ for $n=1,2,3$. Instead of sampling all time instances, we sample partial horizontal slabs which only covers the time steps 1 to 36, 24 to 48, and 30 to 72 for the three areas, respectively. For frontal slab samplings, we sample each area at  different six time steps:  $\mathcal{S}_1^t=\{1, 13, 25, 37, 49, 61\}$, $\mathcal{S}_2^t=\{5, 17, 29, 41, 53, 65\}$, and $\mathcal{S}_3^t=\{9, 21, 33, 45, 57, 69\}$. In the following paragraph, we will verify that the slab sampling schemes satisfy the conditions in Proposition~\ref{proposition:identifiability} when the rank ($F$) of CPD is 7.

First, for each area, $n=1,2,3$, $|\mathcal{S}_n^p|\geq 2$, $|\mathcal{S}_n^t|\geq 2$, so we have sampled subtensors formed in each area. For horizontal slab samplings, let $K_n$ denote the sampled time steps in area $n$. Then we have that $K_1=36$, $K_2=25$, and $K_3=43$. Let $J=5$ denote the number of measurement types. Then
we have that $\min \big\{ \floor{|\mathcal{S}_n^p|} + \floor{J}, \floor{|\mathcal{S}_n^{p}|} + \floor{K_n}, \floor{J} + \floor{K_n}, \log_2(4J|\mathcal{S}_n^t|)\big\} \geq \log_2(4 F)$ for each area $n$, which guarantees that the 
sampled horizontal subtensor in each area admits a unique CPD \cite{Kanatsoulis2020}. The frontal slab samplings in each area guarantee that $\mathcal{P}_n=\mathcal{I}_n$ for any $n=1,2,3$ and  $\bigcup_{n=1}^3 \mathcal{M}_n=\{1,\ldots,J\}$; and  the horizontal slab samplings in each area guarantee that $\bigcup_{n=1}^3 \mathcal{T}_n=\{1,\ldots, K\}$; thus, condition 1) in Proposition~\ref{proposition:identifiability} holds. Next, we check condition 2) in Proposition~\ref{proposition:identifiability}. From the horizontal slab sampling in each area, we have $|\mathcal{M}_m\cap \mathcal{M}_n|=5$ for any area $m$ and area $n$, and $|\mathcal{T}_1\cap \mathcal{T}_2|=13, |\mathcal{T}_1\cap \mathcal{T}_3|=7, |\mathcal{T}_2\cap \mathcal{T}_3|=7$, which implies that condition 2) in   Proposition~\ref{proposition:identifiability} also holds.
\begin{table}
	\renewcommand{\arraystretch}{1}
	\caption{Comparison of the proposed algorithm and centralized solver for slab sampling }\vspace{-10pt}
\label{tab:comparisonofadmmsdf}
	\begin{center}
		\begin{tabular}{l c c c c } 
			\toprule
	            & MAPE($|V|$)  & MAE($\theta$) & MAE($P$) & MAE($Q$) \\
			\midrule
			ADMM &  $0.5922\%$ & $0.7758$ & $0.0279$ & $ 0.0105$ \\
		SDF-NLS & $0.6394\%$ & $0.6457 $ & $0.0314$ & $0.0147$ \\
			Data vs. noisy data & $0.8199\%$ & $0.4756$ & $0.0079$ & $0.0080$  \\
			\bottomrule
		\end{tabular}\vspace{-15pt}
	\end{center}
\end{table}

Table~\ref{tab:comparisonofadmmsdf} shows the performance of our proposed algorithm and the centralized algorithm. The performance results in terms of MAPE and MAE  are comparable for the proposed algorithm and for the centralized solver. The runtime for the proposed algorithm is 16 seconds, whereas that for the centralized solver is
41 seconds, which demonstrates the time efficiency of the proposed algorithm over the centralized solver for the similar performance. \rev{
This can be explained by that the computational complexity for the centralized method in \cite{Zamzam2020} per iteration is ${O}({|\Omega|F})$, and that for our proposed distributed method per iteration is ${O}({{|\Omega|}F/N})$, where $|\Omega|$  and $N$ are the number of measurements in the whole network and the number of partitions, respectively.} The MAPE and MAEs between the true data and the data with 1\% Gaussian noise are also shown in Table~\ref{tab:comparisonofadmmsdf}, which demonstrates the efficacy of the proposed approach in rejecting noise and estimating network states under noisy measurements.

\begin{figure}
    \centering
   \includegraphics[scale=0.5,trim= 80 230 50 240, clip]{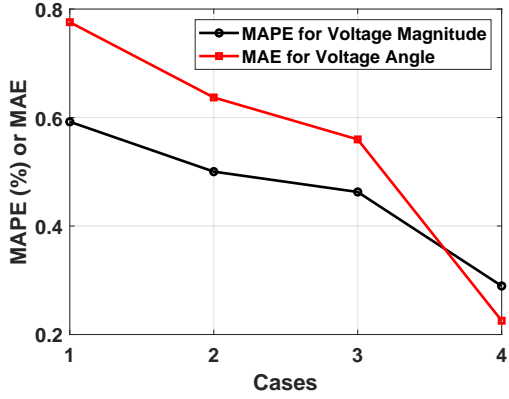}
    \caption{MAPE ($|V|$) and MAE ($\theta$) vs. cases}
     \label{fig:slabsampling1}
 \end{figure}

\begin{figure}
     \centering
   \includegraphics[scale=0.5,trim= 80 240 50 250, clip]{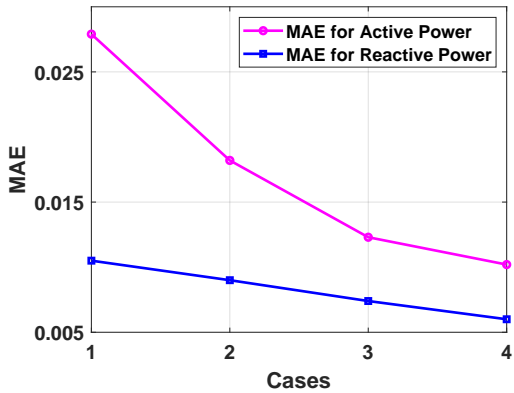}
     \caption{MAE ($P$) and MAE ($Q$) vs. cases}
    \label{fig:slabsampling2}
 \end{figure}

Figs.~\ref{fig:slabsampling1} and \ref{fig:slabsampling2} show that as more frontal slabs are sampled in each area, the MAPE and MAEs decrease. We refer to the aforementioned scenario as Case 1, and we consider the following three cases\footnote{We use \MATLAB notations to denote the set of sampled slabs indices.}:
for Case 2,  $\mathcal{S}_1^t=\{1:6:72\}$, $\mathcal{S}_2^t=\{5:6:72\}$, and $\mathcal{S}_3^t=\{9:6:72\}$; for Case 3, 
$\mathcal{S}_1^t=\{1:3:72\}$, $\mathcal{S}_2^t=\{5:3:72\}$, and $\mathcal{S}_3^t=\{9:3:72\}$; and for Case 4, 
$\mathcal{S}_1^t=\{1:72\}$, $\mathcal{S}_2^t=\{5:72\}$, and $\mathcal{S}_3^t=\{9:72\}$.
Note that the set of sampled time steps for case 4 includes the set of sampled time steps for case 3. Similarly, the set of sampled time steps in case 3  includes the corresponding set for case 2, and so on. For all the cases, the horizontal slabs sampled are the same.


\subsection{Fiber Sampling}
\label{sec:fibersampling}
\begin{figure}
    \centering
  \includegraphics[scale=0.5]{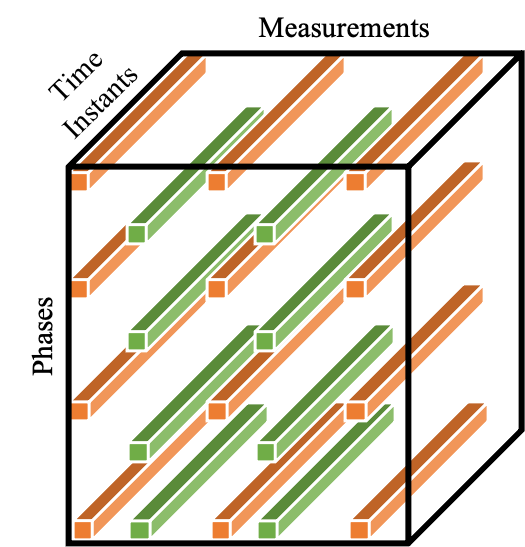}
    \caption{Two-pattern fiber sampling example.}
    \label{fig:fiber}
\end{figure}

For fiber sampling, we use  $\Ten{\bar{X}}_n(\mathcal{S}_n^p,\mathcal{S}_n^m,:)$ to denote the sampled subtensor in area $n$, where $\mathcal{S}_n^p\subseteq\mathcal{I}_n$ denotes the set of sampled phases, and $\mathcal{S}_n^m\subseteq\{1,\ldots, J\}$ denotes the set of  sampled measurement types. 
Let $\mathcal{S}_n^p=\mathcal{I}_n$ for $n=1,2,3$, $\mathcal{S}_1^m =\{1,2,3\}$, $\mathcal{S}_2^m =\{3,4,5\}$, and $\mathcal{S}_3^m =\{3,4\}$; Fig. \ref{fig:fiber} shows an example of a fiber sampling scheme with two patterns. Instead of sampling all the time steps, we sample
the time steps $\{1:2:72\}$, $\{1:72\}$, and $\{1:3:72\}$  for the three areas, respectively. The sampling schemes mean that for all the phases in Area 1, the real voltage, imaginary voltage, and voltage magnitude are sampled  at time steps $\{1,3,5,\ldots, 71\}$; for all the phases  in Area 2, the real voltage, active power, and reactive power are sampled at  time steps $\{1,2,3,\ldots, 72\}$; for all the phases  in Area 3, the voltage magnitude and active power are sampled at  time steps $\{1,4,7,\ldots, 70\}$. Next, we will verify that the fiber sampling schemes satisfy the conditions in Proposition~\ref{proposition:identifiability}  when the rank ($F$) of CPD is 7.

First, for each area $n=1,2,3$, $|\mathcal{S}_n^p|\geq 2$, $|\mathcal{S}_n^m|\geq 2$, so we have sampled subtensors formed in each area.
Let $K_n$ denote the sampled time steps in area $n$. Then we have that $K_1=36$, $K_2=24$, and $K_3=18$. The inequality $\min \big\{ \floor{|\mathcal{S}_n^{p}|} + \floor{|\mathcal{S}_n^{m}|}, \floor{|\mathcal{S}_n^{p}|} + \floor{K}, \floor{|\mathcal{S}_n^{m}|} + \floor{K_n}\big\} \geq \log_2(4 F)$ for $n=1,2,3$ guarantees that the 
sampled  subtensor in each area admits a unique CPD \cite{Kanatsoulis2020}. 
 The sampled phases satisfying $\mathcal{S}_n^p=\mathcal{I}_n$ for $n=1,2,3$ guarantee that $\mathcal{P}_n=\mathcal{I}_n$ for $n=1,2,3$, and the sampled measurement types satisfying $\mathcal{S}_1^m =\{1,2,3\}$, $\mathcal{S}_2^m =\{3,4,5\}$, and $\mathcal{S}_3^m =\{3,4\}$ guarantee that $\bigcup_{n=1}^3 \mathcal{M}_n=\{1,\ldots,J\}$; and  the sampled time steps $\{1:2:72\}, \{1:72\}, \{1:3:72\}$  guarantee that $\bigcup_{n=1}^3 \mathcal{T}_n=\{1,\ldots, K\}$. Thus, condition 1) in Proposition~\ref{proposition:identifiability} holds. Next we check condition 2) in Proposition~\ref{proposition:identifiability}. By that $\mathcal{S}_1^m =\{1,2,3\}$, $\mathcal{S}_2^m =\{3,4,5\}$, $\mathcal{S}_3^m =\{3,4,\}$ , we have that $|\mathcal{M}_1\cap \mathcal{M}_2|=1, |\mathcal{M}_1\cap \mathcal{M}_3|=1,$ and $|\mathcal{M}_2\cap \mathcal{M}_3|=2$. The sampled time steps $\{1:2:72\}, \{1:72\}, \{1:3:72\}$ in the three areas  imply that $|\mathcal{T}_1\cap \mathcal{T}_2|=36, |\mathcal{T}_1\cap \mathcal{T}_3|=12, |\mathcal{T}_2\cap \mathcal{T}_3|=24$, so condition 2) in   Proposition~\ref{proposition:identifiability} also holds. Thus, all the conditions in Proposition~\ref{proposition:identifiability} are satisfied.

\begin{table}
	\renewcommand{\arraystretch}{1}
	\caption{Comparison of the proposed algorithm and the centralized solver for fiber sampling}\vspace{-10pt}
\label{tab:comparisonofadmmsdffiber}
	\begin{center}
		\begin{tabular}{l c c c c } 
			\toprule
	            & MAPE($|V|$)  & MAE($\theta$) & MAE($P$) & MAE($Q$) \\
			\midrule
			ADMM &  $0.3083\%$ & $0.9576$ & $0.0547$ & $ 0.0157$ \\
		SDF-NLS & $0.2412\%$ & $0.9734 $ & $0.1582$ & $0.0421$\\
		\bottomrule
		\end{tabular}\vspace{-15pt}
	\end{center}
\end{table}

Table~\ref{tab:comparisonofadmmsdffiber} shows that the MAPE of the voltage magnitude for the centralized solver is slightly better than the proposed algorithm, and the MAE of the voltage angle is comparable, but the MAEs of the active power and the reactive power for the proposed algorithm are much better than those for the centralized solver. The runtime for the the proposed algorithm is 15 seconds, and that for the centralized solver is 41 seconds. The results show that the proposed algorithm slightly outperforms the centralized solver in terms of both accuracy and time efficiency while reducing the computational time and not requiring the data to be shared from all controllers.

\rev{It is worth noting that for the fiber sampling, we only sample three phases at the slack bus with real voltage, imaginary voltage, and voltage magnitude known, and the other phases in the first area with voltage magnitude, active power, and reactive power known, for the other two areas, the assumption is the same with the manuscript, and we get MAPE for the voltage magnitude around $0.4023\%$ (a little bit worse than the above scenario), MAE for the voltage angle around $1.0484$ (a little bit worse than the above scenario), and MAEs for the active power and reactive power are $0.0534$ and $0.0164$ (comparable with the above scenario). }

\subsection{Mixed Sampling}
In this section, we demonstrate that our proposed algorithm works for  mixed sampling schemes and that our proposition can be used to check mixed sampling schemes too. We assume that Area 1 has horizontal and frontal slab sampling schemes, and areas 2 and 3 have fiber sampling schemes. We use the same notations with sections~\ref{sec:slabsampling} and \ref{sec:fibersampling}. For horizontal slab sampling, we assume that $|\mathcal{S}_1^p|=8$, and  we sample only at time steps $\{1:9:72\}=\{1,10,19,28,37,46,55,64\}$.
For frontal slab samplings, we sample   at time steps $18,36,56,72$. For fiber sampling, we assume that $\mathcal{S}_n^p=\mathcal{I}_n$ for $n=2,3$, $\mathcal{S}_2^m =\{3,4,5\}$, and $\mathcal{S}_3^m =\{1,2,3\}$. Instead of sampling all the time steps, we sample
the time steps $\{1:72\}$,  and $\{1:3:72\}$  for the two areas, respectively. Similar to sections~\ref{sec:slabsampling} and \ref{sec:fibersampling}, we can check that the mixed sampling schemes satisfy the conditions in Proposition~\ref{proposition:identifiability} when the rank ($F$) of CPD is 7.

Table~\ref{tab:comparisonofadmmsdfmixed} shows that the MAPE of the voltage magnitude and the MAEs of the voltage angle  for both the proposed algorithm and the centralized solver are comparable, but the MAEs of the active power and the reactive power for the proposed algorithm are much better than those for the centralized solver. 
\begin{table}
	\renewcommand{\arraystretch}{1}
	\caption{Comparison of the proposed algorithm and the centralized solver for mixed sampling }\vspace{-10pt}
\label{tab:comparisonofadmmsdfmixed}
	\begin{center}
		\begin{tabular}{l c c c c } 
			\toprule
	            & MAPE($|V|$)  & MAE($\theta$) & MAE($P$) & MAE($Q$) \\
			\midrule
			ADMM &  $0.5054\%$ & $0.9318$ & $0.0326$ & $ 0.0168$ \\
		SDF-NLS & $0.6431\%$ & $0.9290 $ & $0.1459$ & $0.0382$\\
		\bottomrule
		\end{tabular}\vspace{-15pt}
	\end{center}
\end{table}

\section{Concluding Remarks and Future Research}
\label{sec:conclusion}
This paper considered a  model-free tensor completion method using CPD for DSSE. By dividing the network into multiple areas, the measured physical quantities at buses located in the same area were processed by an area controller. A three-way tensor was constructed to collect the measured quantities analyzed locally to recover the full state information of the network. A distributed closed-form ADMM algorithm was developed to obtain the low-rank factors of the whole network state tensor where information exchange happened only between neighboring areas.  The convergence properties of
the distributed algorithm and the sufficient conditions on the number of samples for each smaller network that guarantee the identifiability of the factors of the state tensor were investigated.
In  addition,  simulations  were  carried out 
on the  IEEE 123-bus  distribution  feeder to  verify the identifiability conditions and to demonstrate the efficacy of the distributed ADMM.   


%





\ifCLASSOPTIONcaptionsoff
  \newpage
\fi




\appendices
\rev{
\section{Proof of Proposition 1}
To prove this proposition, we first invoke two results of the tensor algebra literature. The first theorem introduced general conditions under which the CPD decomposition of a tensor is essentially unique almost surely. 
\begin{theorem}\cite{Chiantini2012}
Let $\Ten{X} = {\cpd{\A}{\B}{\C}}$ with ${\bf A}\in \mathbb{R}^{I\times F}$, ${\bf B}\in \mathbb{R}^{J\times F}$, and ${\bf C}\in \mathbb{R}^{K\times F}$. Assume that ${\bf A}$, ${\bf B}$, and ${\bf C}$ are drawn from some joint absolutely continuous distribution. Also, assume $I \geq J \geq K$ without loss of generality. If $F\leq 2^{\lfloor \log_2 J\rfloor + \lfloor \log_2 K\rfloor -2}$, then the decomposition of $\Ten{X}$ in terms of ${\bf A}$, ${\bf B}$, and ${\bf C}$ is essentially unique almost surely.
\end{theorem}
Essential uniqueness means that the decomposition is unique up to a common rotation and scaling and counter-scaling. In other words, any decomposition of $\Ten{X}$ is a rotated and scaled and counter-scaled version of the unique factors, i.e.,
\begin{align}
    \overline{\bf A} = {\bf A} \boldsymbol{\Pi} \boldsymbol{\Gamma}_1,\quad
    \overline{\bf B} = {\bf B} \boldsymbol{\Pi} \boldsymbol{\Gamma}_2,\quad 
    \overline{\bf C} = {\bf C} \boldsymbol{\Pi} \boldsymbol{\Gamma}_3,
\end{align}
where $\boldsymbol{\Pi}$ denotes a common rotation matrix, and $\boldsymbol{\Gamma}_i$ are diagonal matrices satisfying $\boldsymbol{\Gamma}_1 \boldsymbol{\Gamma}_2 \boldsymbol{\Gamma}_3 = \boldsymbol{I}$. Then, we introduce the second result from the literature which ensures that the uniqueness property is preserved even if the underlying factors are multiplied by a deterministic matrix.
\begin{theorem}\cite{kanatsoulis2018hyperspectral}
Let $\tilde{\bf Z} = {\bf P} {\bf Z} $, where the elements of ${\bf Z}$ are drawn from an absolutely continuous joint distribution with respect to the Lebesgue measure in $\mathbb{R}^{IF}$, and ${\bf P}$ is a deterministic matrix of size $I' \times I$ which has a full row rank. Then, the joint distribution of the elements of $\tilde{\bf Z}$ is absolutely continuous with respect to the Lebesgue measure in $\mathbb{R}^{I'F}$.
\end{theorem}

This result guarantees that the essential uniqueness property is preserved if the original tensor is distributed as we propose in this paper. Before delving into the proof of our proposition, we present an example of the recovery process of the original factors in a simple scenario.

Assume the original tensor $\Ten{X}\in\mathbb{R}^{10\times 4\times 4}$ is divided into two partitions. The first partition collects all measurements for the first five nodes. The second partition collects all measurements of the last five nodes of the network. After performing decomposing the two partitions, the resulting factors would be
\begin{align}\label{eq:factors-a}
    &\overline{\bf A}^{(1)} = {\bf P}_a^{(1)} {\bf A} \boldsymbol{\Pi}^{(1)} \boldsymbol{\Gamma}_a^{(1)}, \quad 
    \overline{\bf A}^{(2)} = {\bf P}_a^{(2)} {\bf A} \boldsymbol{\Pi}^{(2)} \boldsymbol{\Gamma}_a^{(2)},\\\label{eq:factors-b}
    &\overline{\bf B}^{(1)} = {\bf B} \boldsymbol{\Pi}^{(1)} \boldsymbol{\Gamma}_b^{(1)}, \quad 
    \qquad\overline{\bf B}^{(2)} = {\bf B} \boldsymbol{\Pi}^{(2)} \boldsymbol{\Gamma}_b^{(2)}, \\\label{eq:factors-c}
    &\overline{\bf C}^{(1)} = {\bf C} \boldsymbol{\Pi}^{(1)} \boldsymbol{\Gamma}_c^{(1)}, \quad 
    \qquad\overline{\bf C}^{(2)} = {\bf C} \boldsymbol{\Pi}^{(2)} \boldsymbol{\Gamma}_c^{(2)},
\end{align}
where ${\bf P}_a^{(1)}:= [{\bf I}_5\quad {\bf 0}]$ and ${\bf P}_a^{(2)}:= [{\bf 0}\quad {\bf I}_5]$. Let us denote the inverse of $\boldsymbol{\Gamma}_{b}^{(i)}$ by $\boldsymbol{\Lambda}_b^{(i)}$, similar for $\boldsymbol{\Gamma}_c^{(i)}$. Then, to reconstruct ${\bf A}$ from the copies, we need to resolve the rotations $\boldsymbol{\Pi}^{(1)}$ and $\boldsymbol{\Pi}^{(2)}$ as well as the scaling ambiguities $\boldsymbol{\Gamma}_a^{(1)}$ and $\boldsymbol{\Gamma}_a^{(2)}$. We aim at recovering 
\begin{align}\notag
    {\bf A} = \begin{bmatrix} {\bf P}_a^{(1)} {\bf A}\\
    {\bf P}_a^{(2)} {\bf A}
    \end{bmatrix} = \begin{bmatrix} \overline{\bf A}^{(1)} \boldsymbol{\Lambda}_a^{(1)} {\boldsymbol{\Pi}^{(1)}}^T\\
    \overline{\bf A}^{(2)} \boldsymbol{\Lambda}_a^{(2)} {\boldsymbol{\Pi}^{(2)}}^T
    \end{bmatrix}= \begin{bmatrix} \overline{\bf A}^{(1)} \boldsymbol{\Gamma}_c^{(1)} \boldsymbol{\Gamma}_b^{(1)} {\boldsymbol{\Pi}^{(1)}}^T\\
    \overline{\bf A}^{(2)} \boldsymbol{\Gamma}_c^{(2)} \boldsymbol{\Gamma}_b^{(2)} {\boldsymbol{\Pi}^{(2)}}^T
    \end{bmatrix}.
\end{align}
If we matricize the first partition of the tensor, we get
\begin{align}\nonumber
    \Ten{X}_{1} &= (\overline{\bf C}^{(1)} \odot \overline{\bf B}^{(1)}) {\overline{\bf A}^{(1)}}^T\\\nonumber
    &= (\overline{C}^{(2)} \boldsymbol{\Lambda}_c^{(2)} {\boldsymbol{\Pi}^{(2)}}^T \boldsymbol{\Pi}^{(1)} \boldsymbol{\Gamma}_c^{(1)} \\\nonumber &\qquad\odot \overline{B}^{(2)} \boldsymbol{\Lambda}_b^{(2)} {\boldsymbol{\Pi}^{(2)}}^T \boldsymbol{\Pi}^{(1)} \boldsymbol{\Gamma}_b^{(1)}) {\overline{A}^{(1)}}^T\\\nonumber
    &=(\overline{\bf C}^{(2)} \boldsymbol{\Lambda}_c^{(2)} \tilde{\boldsymbol{\Gamma}}_c^{(1)}{\boldsymbol{\Pi}^{(2)}}^T \boldsymbol{\Pi}^{(1)} \\\nonumber &\qquad \odot \overline{\bf B}^{(2)} \boldsymbol{\Lambda}_b^{(2)} \tilde{\boldsymbol{\Gamma}}_b^{(1)}{\boldsymbol{\Pi}^{(2)}}^T \boldsymbol{\Pi}^{(1)}) {\overline{\bf A}^{(1)}}^T\\\nonumber
    &=(\overline{\bf C}^{(2)} \odot \overline{\bf B}^{(2)}) \boldsymbol{\Lambda}_b^{(2)} \tilde{\boldsymbol{\Gamma}}_b^{(1)} \boldsymbol{\Lambda}_c^{(2)} \tilde{\boldsymbol{\Gamma}}_c^{(1)} {\boldsymbol{\Pi}^{(2)}}^T \boldsymbol{\Pi}^{(1)} {\overline{\bf A}^{(1)}}^T\\\nonumber
    &=(\overline{\bf C}^{(2)} \odot \overline{\bf B}^{(2)}) \boldsymbol{\Lambda}_b^{(2)}  \boldsymbol{\Lambda}_c^{(2)}  {\boldsymbol{\Pi}^{(2)}}^T \boldsymbol{\Pi}^{(1)} {\boldsymbol{\Gamma}}_b^{(1)}{\boldsymbol{\Gamma}}_c^{(1)}{\overline{\bf A}^{(1)}}^T
\end{align}
where in the first step we do simple manipulations of \eqref{eq:factors-b} and \eqref{eq:factors-c}. Then, in the second and fifth steps, we utilize the property that ${\boldsymbol{\Pi}^{(2)}}^T \boldsymbol{\Pi}^{(1)} \boldsymbol{\Gamma}_c^{(1)} = \tilde{\boldsymbol{\Gamma}}_c^{(1)} {\boldsymbol{\Pi}^{(2)}}^T \boldsymbol{\Pi}^{(1)} $ where $\tilde{\boldsymbol{\Gamma}}_c^{(1)}$ is a diagonal matrix that has the same elements with ${\boldsymbol{\Gamma}}_c^{(1)}$ but rotated according to $ {\boldsymbol{\Pi}^{(2)}}^T \boldsymbol{\Pi}^{(1)} $. Then, the fourth step utilizes Khatri-Rao product properties. Using the results from \cite{jiang2001almost}, the matrix $(\overline{\bf C}^{(2)} \odot \overline{\bf B}^{(2)})$ is full column rank almost surely. Therefore, from the matricization of the first tensor, we can recover $\tilde{\bf A}:=\overline{\bf A}^{(1)} \boldsymbol{\Gamma}_c^{(1)} \boldsymbol{\Gamma}_b^{(1)} {\boldsymbol{\Pi}^{(1)}}^T \boldsymbol{\Pi}^{(2)} \boldsymbol{\Lambda}_b^{(2)} \boldsymbol{\Lambda}_c^{(c)}$.

Using the recovered $\tilde{\bf A}$, we can write 
\begin{align}
    {\bf A} \boldsymbol{\Pi}^{(2)} \boldsymbol{\Lambda}_b^{(2)} \boldsymbol{\Lambda}_c^{(c)} = \begin{bmatrix} \tilde{\bf A}\\
    \overline{\bf A}^{(2)}
    \end{bmatrix}
\end{align}
This represents a rotated and scaled version of the generating factor ${\bf A}$. This way we combine the identified rows from $\overline{\bf A}^{(1)}$ and $\overline{\bf A}^{(2)}$ which have different rotations and scaling ambiguities. In general, the proof of the theorem details the conditions under which the previous procedure can be applied to recover the factors of the tensor up to a common permutation and scaling and counter scaling ambiguities.
\begin{proof}
The condition that the sampled subtensor $\Ten{\bar{X}}_n$ for every area $n$ admits a unique CPD guarantees that the factors of the subtensor constructed from the sampled subtensor at each area can be identified, up to column permutation and scaling. 
In order to recover the factors of the original tensor, the low-rank factors identified from each area need to be combined. Let the identified factors from area $n$ be $(\check{\A}_n, \check{\B}_n, \check{\C}_n)$. In this case, we relate these factors to the original factors as follows:
\begin{align}
    \check{\Ten{X}_n} \approx \big[[ {\bf P}_{n, a} \boldsymbol{\Pi}_n {\bf S}_{n, a} \A, {\bf P}_{n, b} \boldsymbol{\Pi}_n {\bf S}_{n, b} \B, {\bf P}_{n, c} \boldsymbol{\Pi}_n {\bf S}_{n, c} \C]\big]
\end{align}
where ${\bf P}_{n, a}$, ${\bf P}_{n, b}$, and ${\bf P}_{n, c}$ represent the scaling diagonal matrices satisfying ${\bf P}_{n, a} {\bf P}_{n, b} {\bf P}_{n, c} = \boldsymbol{I}$. Also, $\boldsymbol{\Pi}_n$ denotes the common rotation of factors and the selected rows of $\A$, $\B$, and $\C$ are sampled using the sampling matrices ${\bf S}_{n, a}$, ${\bf S}_{n, b}$, and ${\bf S}_{n, c}$, respectively. Assume that two areas $m$ and $n$ are mutually identifiable. Then, in order to combine the rows available in $(\check{\A}_n, \check{\B}_n, \check{\C}_n)$ and $(\check{\A}_m, \check{\B}_m, \check{\C}_m)$, we need to recover the rotation $\boldsymbol{\Pi}_n$ and $\boldsymbol{\Pi}_m$ as well as the scaling ambiguities. We consider the two possible cases:\\
i) $|\mathcal{M}_n\cap \mathcal{M}_m|\geq 2 \ \text{and}\  |\mathcal{T}_n\cap \mathcal{T}_m|\geq 1$. In this case, we select any two common rows of $\mathcal{M}_n$ and $\mathcal{M}_m$. By calculating the element-wise ratios between the two rows on $\check{\B}_n$ and $\check{B}_c$, we can resolve the permutation ambiguity between the two sets of factors. We apply the counter permutations on all $\A$ and $\C$ factors as well. We can then absorb the scaling in $\A$ by setting the elements of the one of the common rows of $\check{\B}_n$ and $\check{\B}_m$ as well as the common row of $\check{\C}_n$ and $\check{\C}_m$ to be the same. Then, we proceed by just combining the factors since the sampling matrices are known.\\
ii) $|\mathcal{T}_n\cap \mathcal{T}_m|\geq 2 \ \text{and}\  |\mathcal{M}_n\cap \mathcal{M}_m|\geq 1$. In this case, we can combine the factors similar to the previous case but using the common rows of $\check{\C}_n$ and $\check{\C}_m$ instead of $\check{\B}_n$ and $\check{\B}_m$ and vice versa.

Condition 1) guarantees that all the rows of $\A_n$, $\B$, and $\C$ can be identified from all partitions combined. Condition 2) implies that either i) $|\mathcal{M}_n\cap \mathcal{M}_m|\geq 2 \ \text{and}\  |\mathcal{T}_n\cap \mathcal{T}_m|\geq 1$ hold, or ii) $|\mathcal{T}_n\cap \mathcal{T}_m|\geq 2 \ \text{and}\  |\mathcal{M}_n\cap \mathcal{M}_m|\geq 1$ hold for area $n$ and at least one area $m$. Then, condition 2) also implies that all areas are connected, and hence all  rows of $\A$, $\B$, and $\C$ can be recovered. Therefore, $\Ten{X}$ is recovered, which completes the proof.
\end{proof}
}

\section{Proof of Theorem 1}
The proof follows the steps of the proof of Theorem 2.1  in \cite{Xu2012}. First, we prove that $\mathcal{L}(\{\A_i\}, \{\B_i\}, \{\C_i\}, \{\B_{ij}\}, \{\C_{ij}\}, \{{\bf{\Gamma}}_{ij}\}, \{{\bf{\Lambda}}_{ij}\})$ is strongly convex with respect to each variable of 
${\bf a}_i, {\bf b}_i, {\bf c}_i,\B_{ij}$, and $\C_{ij}$. 
Because 
\begin{align*}
  &\|\Ten{M}_i * \big(\Ten{X}_i - {\cpd{\A_i}{\B_i}{\C_i}} \big)\|_F^2\\
    &=\|\underline{\bf w}_i^{(1)} * \underline{\bf x}_i^{(1)}-\text{diag}(\underline{\bf w}_i^{(1)})((\C_i\odot \B_i)\otimes {\bf I}_{_{I_i\times I_i}}){\bf a}_i\|_2^2\\
    &=\|\underline{\bf w}_i^{(2)} * \underline{\bf x}_i^{(2)}-\text{diag}(\underline{\bf w}_i^{(2)})((\C_i\odot \A_i)\otimes {\bf I}_{_{J\times J}}){\bf b}_i\|_2^2\\
    &=\|\underline{\bf w}_i^{(3)} * \underline{\bf x}_i^{(3)}-\text{diag}(\underline{\bf w}_i^{(3)})((\B_i\odot \A_i)\otimes {\bf I}_{_{K\times K}}){\bf c}_i\|_2^2,
\end{align*}
and 
\begin{align*}
    \|\B_{i}-\B_{ij}+{\bf{\Gamma}}_{ij}\|_F^2&=\|{\bf b}_i-\text{vec}(\B_{ij}-{\bf{\Gamma}}_{ij})\|_2^2\\
     \|\C_{i}-\C_{ij}+{\bf{\Lambda}}_{ij}\|_F^2&=\|{\bf c}_i-\text{vec}(\C_{ij}-{\bf{\Lambda}}_{ij})\|_2^2,
\end{align*}
which show that $\mathcal{L}$ is a strongly convex quadratic function with respect to each variable of 
${\bf a}_i, {\bf b}_i, {\bf c}_i,\B_{ij}$, and $\C_{ij}$, respectively. \revision{By the form of $\mathcal{L}$ and the assumption in the theorem, we have that $\mathcal{L}$ is bounded. By the fact that  $\mathcal{L}$ is bounded and that $\mathcal{L}$ is  strongly convex  to each variable of 
${\bf a}_i, {\bf b}_i, {\bf c}_i,\B_{ij}$, and $\C_{ij}$, respectively, we have the following inequalities hold:
\begin{subequations}
\label{ineq:proof1}
\begin{equation}
    \mathcal{L}({\bf a}_i^{k})- \mathcal{L}({\bf a}_i^{k+1})\geq \|{\bf a}_i^{k+1}-{\bf a}_i^{k}\|_2^2,
\end{equation}
\begin{equation}
    \mathcal{L}({\bf b}_i^{k})- \mathcal{L}({\bf b}_i^{k+1})\geq \mu\|{\bf b}_i^{k+1}-{\bf b}_i^{k}\|_2^2,
\end{equation}
\begin{equation}
    \mathcal{L}({\bf c}_i^{k})- \mathcal{L}({\bf c}_i^{k+1})\geq \lambda\|{\bf c}_i^{k+1}-{\bf c}_i^{k}\|_2^2,
\end{equation}
\begin{equation}
    \mathcal{L}(\B_{ij}^{k})- \mathcal{L}(\B_{ij}^{k
    +1})\geq \mu\|\B_{ij}^{k+1}-\B_{ij}^{k}\|_F^2,
\end{equation}
\begin{equation}
    \mathcal{L}(\C_{ij}^{k})- \mathcal{L}(\C_{ij}^{k+1})\geq \lambda\|\C_{ij}^{k+1}-\C_{ij}^{k}\|_F^2.
\end{equation}
\end{subequations}

Let $a := \min\{1, \mu, \lambda\}$, $\bZ^k = (\{{\bf a}_i^k\}, \{{\bf b}_i^k\}, \{{\bf c}_i^k\}, \{\B_{ij}^k\}, \{\C_{ij}^k\}, \{{\bf{\Gamma}}_{ij}^k\}, \{{\bf{\Lambda}}_{ij}^k\})$, and $\bZ_1^k = (\{{\bf a}_i^k\}, \{{\bf b}_i^k\}, \{{\bf c}_i^k\}, \{\B_{ij}^k\}, \{\C_{ij}^k\})$.  By (\ref{ineq:proof1}), we have that 
\begin{align*}
    &\mathcal{L}(\bZ^k) - \mathcal{L}(\bZ^{k+1})\\
    =&\mathcal{L}(\bZ_1^k, \{{\bf{\Gamma}}_{ij}^k\}, \{{\bf{\Lambda}}_{ij}^k\})-\mathcal{L}(\bZ_1^{k+1}, \{{\bf{\Gamma}}_{ij}^{k+1}\}, \{{\bf{\Lambda}}_{ij}^{k+1}\})\\
    =&\mathcal{L}(\bZ_1^k, \{{\bf{\Gamma}}_{ij}^k\}, \{{\bf{\Lambda}}_{ij}^k\})-
    \mathcal{L}(\bZ_1^{k+1}, \{{\bf{\Gamma}}_{ij}^{k}\}, \{{\bf{\Lambda}}_{ij}^{k}\})+\\
    &\mathcal{L}(\bZ_1^{k+1}, \{{\bf{\Gamma}}_{ij}^{k}\}, \{{\bf{\Lambda}}_{ij}^{k}\})-\mathcal{L}(\bZ_1^{k+1}, \{{\bf{\Gamma}}_{ij}^{k+1}\}, \{{\bf{\Lambda}}_{ij}^{k+1}\})\\
    \geq & a\sum\limits_{i=1}^N(\|{\bf a}_i^{k+1}-{\bf a}_i^{k}\|_2^2+\|{\bf b}_i^{k+1}-{\bf b}_i^{k}\|_2^2+\|{\bf c}_i^{k+1}-{\bf c}_i^{k}\|_2^2+\\
    &\|\B_{ij}^{k+1}-\B_{ij}^{k}\|_F^2+\|\C_{ij}^{k+1}-\C_{ij}^{k}\|_F^2)-\max\{\mu, \lambda\}\\
    &\sum\limits_{i=1}^N\sum\limits_{j\in\mathcal{N}(i)}(\|{\bf{\Gamma}}_{ij}^{k+1}-{\bf{\Gamma}}_{ij}^{k}\|_F^2+\|{\bf{\Lambda}}_{ij}^{k+1}-{\bf{\Lambda}}_{ij}^{k}\|_F^2).
\end{align*}
Taking summation of the above inequality for $k$ from 0 to $\infty$ and recalling that $\mathcal{L}$ is bounded, we have that  
\begin{align*}
    &a\sum\limits_{k=0}^{\infty}\sum\limits_{i=1}^N(\|{\bf a}_i^{k+1}-{\bf a}_i^{k}\|_2^2+\|{\bf b}_i^{k+1}-{\bf b}_i^{k}\|_2^2+\|{\bf c}_i^{k+1}-{\bf c}_i^{k}\|_2^2+\\
    &\|\B_{ij}^{k+1}-\B_{ij}^{k}\|_F^2+\|\C_{ij}^{k+1}-\C_{ij}^{k}\|_F^2)-\max\{\mu, \lambda\}\\
    &\sum\limits_{k=0}^{\infty}\sum\limits_{i=1}^N\sum\limits_{j\in\mathcal{N}(i)}(\|{\bf{\Gamma}}_{ij}^{k+1}-{\bf{\Gamma}}_{ij}^{k}\|_F^2+\|{\bf{\Lambda}}_{ij}^{k+1}-{\bf{\Lambda}}_{ij}^{k}\|_F^2)<\infty.
\end{align*}
By (\ref{ineq:boundassumption}), we have that 
\begin{align*}
    &\sum\limits_{k=0}^{\infty}\sum\limits_{i=1}^N(\|{\bf a}_i^{k+1}-{\bf a}_i^{k}\|_2^2+\|{\bf b}_i^{k+1}-{\bf b}_i^{k}\|_2^2+\|{\bf c}_i^{k+1}-{\bf c}_i^{k}\|_2^2+\\
    &\|\B_{ij}^{k+1}-\B_{ij}^{k}\|_F^2+\|\C_{ij}^{k+1}-\C_{ij}^{k}\|_F^2)<\infty,
\end{align*}
which implies that (\ref{eq:primalconvergence1}) holds for $i=1,\cdots, N$. And for any $i=1,\cdots, N, j\in\mathcal{N}(i)$, by  (\ref{ineq:boundassumption}) we directly have (\ref{eq:primalconvergence2}) hold.
} 
\begin{subequations}
\label{eq:convergencecond}
\begin{equation}
\label{eq:primalconvergence1}
    ({\bf a}_i^{k+1}, {\bf b}_i^{k+1}, {\bf c}_i^{k+1})- ({\bf a}_i^{k}, {\bf b}_i^{k}, {\bf c}_i^{k})\longrightarrow 0,
\end{equation}
\begin{equation}
\label{eq:primalconvergence2}
     ({\bf{\Gamma}}_{ij}^{k+1}, {\bf{\Lambda}}_{ij}^{k+1})-  ({\bf{\Gamma}}_{ij}^{k}, {\bf{\Lambda}}_{ij}^{k})\longrightarrow 0. 
\end{equation}
\end{subequations}
Then we prove that any accumulation point of $\{\bZ^k\}$ generated by the distributed ADMM for problem 
(\ref{eq:auxidistributedCPDmodel}) is a KKT point satisfying (\ref{eq:KKTconditions}).
By the updates of ${\bf a}_i, {\bf b}_i, {\bf c}_i, {\bf{\Gamma}}_{ij}, {\bf{\Lambda}}_{ij}$, we have:
\begin{subequations}
\label{eq:primaldualupdates}
\begin{equation}
\label{eq:primalupdate1}
    (\D_i^{kT}\D_i^k)({\bf a}_i^{k+1}- {\bf a}_i^{k})
    =\D_i^{kT}(\underline{\bf w}_i^{(1)}* \underline{\bf x}_i^{(1)}-\D_i^k{\bf a}_i^{k}),
\end{equation}
\begin{equation}
\nonumber
    (\E_i^{kT}\E_i^k+\mu\sum_{j\in\mathcal{N}(i)}{\bf I}_{{FJ\times FJ}})({\bf b}_i^{k+1}-{\bf b}_i^{k})
    = 
\end{equation}
\begin{equation}
    \nonumber
    -(\E_i^{kT}\E_i^k+\mu\sum_{j\in\mathcal{N}(i)}{\bf I}_{{FJ\times FJ}}){\bf b}_i^{k}+
    \E_i^{kT} (\underline{\bf w}_i^{(2)}* \underline{\bf x}_i^{(2)})
\end{equation}
\begin{equation}
    \label{eq:primalupdate2}
    +\mu\sum_{j\in\mathcal{N}(i)}\text{vec}(\B_{ij}^{k}-{\bf{\Gamma}}_{ij}^{k}),
\end{equation}
\begin{equation}
\nonumber
    (\bG_i^{kT}\bG_i^k+\mu\sum_{j\in\mathcal{N}(i)}{\bf I}_{{Fk\times Fk}})({\bf c}_i^{k+1}-{\bf c}_i^{k})
    =
\end{equation}
\begin{equation}
    \nonumber
     -(\bG_i^{kT}\bG_i^k+\mu\sum_{j\in\mathcal{N}(i)}{\bf I}_{{Fk\times Fk}}){\bf c}_i^{k}+\bG_i^{kT} (\underline{\bf w}_i^{(3)}* \underline{\bf x}_i^{(3)})
\end{equation}
\begin{equation}
    \label{eq:primalupdate3}
    +\mu\sum_{j\in\mathcal{N}(i)}\text{vec}(\C_{ij}^{k}-{\bf{\Lambda}}_{ij}^{k}),
\end{equation}
\begin{equation}
    \label{eq:dualupdate1}
    {\bf{\Gamma}}_{ij}^{k+1} - {\bf{\Gamma}}_{ij}^k ={(\B_i^{k+1}-\B_j^{k+1})}/{2}, j\in\mathcal{N}(i),
\end{equation}
\begin{equation}
\label{eq:dualupdate2}
{\bf{\Lambda}}_{ij}^{k+1}- {\bf{\Lambda}}_{ij}^k= {(\C_i^{k+1}-\C_j^{k+1})}/{2}, j\in\mathcal{N}(i).
\end{equation}
\end{subequations}
By (\ref{eq:convergencecond}), we have that the left-hand and right-hand sides of (\ref{eq:primaldualupdates}) all go to zero, i.e.:
\begin{subequations}
\label{eq:KKTconvergence}
\begin{equation}
    \label{eq:KKTconvergence1}
    \D_i^{kT}(\underline{\bf w}_i^{(1)}* \underline{\bf x}_i^{(1)}-\D_i^k{\bf a}_i^{k})\longrightarrow 0
\end{equation}
\begin{equation}
    \nonumber
    -(\E_i^{kT}\E_i^k+\mu\sum_{j\in\mathcal{N}(i)}{\bf I}_{{FJ\times FJ}}){\bf b}_i^{k}+
    \E_i^{kT} (\underline{\bf w}_i^{(2)}* \underline{\bf x}_i^{(2)})
\end{equation}
\begin{equation}
    \label{eq:KKTconvergence2}
    +\mu\sum_{j\in\mathcal{N}(i)}\text{vec}(\B_{ij}^{k}-{\bf{\Gamma}}_{ij}^{k})\longrightarrow 0,
\end{equation}
\begin{equation}
    \nonumber
     -(\bG_i^{kT}\bG_i^k+\mu\sum_{j\in\mathcal{N}(i)}{\bf I}_{{Fk\times Fk}}){\bf c}_i^{k}+\bG_i^{kT} (\underline{\bf w}_i^{(3)}* \underline{\bf x}_i^{(3)})
\end{equation}
\begin{equation}
    \label{eq:KKTconvergence3}
    +\mu\sum_{j\in\mathcal{N}(i)}\text{vec}({\C_{ij}^{k}}-{\bf{\Lambda}}_{ij}^{k})\longrightarrow 0,
\end{equation}
\begin{equation}
    \label{eq:KKTconvergence4}
    {\B_i^{k+1}-\B_j^{k+1}}\longrightarrow 0, j\in\mathcal{N}(i),
\end{equation}
\begin{equation}
    \label{eq:KKTconvergence5}
    {\C_i^{k+1}-\C_j^{k+1}}\longrightarrow 0, j\in\mathcal{N}(i).
\end{equation}
\end{subequations}
\rev{Then for any limit point $\hat{\bZ}_1 = (\{\hat{{\bf a}}_i\}, \{\hat{{\bf b}}_i\}, \{\hat{{\bf c}}_i\}, \{\hat{\B}_{ij}\}, \{\hat{\C}_{ij}\})$ of sequence $\{{\bZ}_1^k\}$, there exists a subsequence $\{{\bZ}_1^{n_k}\}$ converging to $\hat{\bZ}_1$. The boundedness of $(\{{\bf{\Gamma}}_{ij}^k\},\{{\bf{\Lambda}}_{ij}^k\})$ implies the existence of a subsequence $(\{{\bf{\Gamma}}_{ij}^{l_k}\},\{{\bf{\Lambda}}_{ij}^{l_k}\})$ of $(\{{\bf{\Gamma}}_{ij}^{n_k}\},\{{\bf{\Lambda}}_{ij}^{n_k}\})$ converging to some point $(\{\hat{{\bf{\Gamma}}}_{ij}\},\{\hat{{\bf{\Lambda}}}_{ij}\})$. Hence, $(\hat{\bZ}_1,\{\hat{{\bf{\Gamma}}}_{ij}\},\{\hat{{\bf{\Lambda}}}_{ij}\})$ is a limit point of $\{\bf{Z}^k\}$. We have that the  five equations in the KKT conditions (\ref{eq:KKTconditions}) are satisfied at the limit point $(\hat{\bZ}_1,\{\hat{{\bf{\Gamma}}}_{ij}\},\{\hat{{\bf{\Lambda}}}_{ij}\})$. This completes the proof.
}

\rev{
\section{Performance on a Real Feeder}
\begin{table}
	\renewcommand{\arraystretch}{1}
	\captionsetup{labelfont={color=black},font={color=black}}
	\caption{Performance on a real utility feeder system}\vspace{-10pt}
\label{tab:realutility}
	\begin{center}
		\begin{tabular}{l c c c c } 
			\toprule
	            & \res{MAPE($|V|$)}  & \res{MAE($\theta$)} & \res{MAE($P$)} & \res{MAE($Q$)} \\
			\midrule
		  &\res{$0.8833\%$} & \res{$0.2552$} & \res{$0.0197$} & \res{$0.0028$} \\
		\bottomrule
		\end{tabular}\vspace{-15pt}
	\end{center}
\end{table}
In order to demonstrate the scalability of our algorithm, we
implement it on a 2576-phase feeder from our utility partner with a ten-area
partition for the following slab sampling scenario: for horizontal slab sampling, we sample ten
phases randomly from nonzero load phases for each area at fifteen time steps; for frontal slab
sampling, we sample each area at different six time steps. Instead of considering the zero active and reactive power injections to be known, we leave them unkown, and the percentage of known measurements is 7.61\%. The performance in terms of MAPE and MAEs are shown in Table \ref{tab:realutility}.

The results in terms of MAPE and
MAEs are similar to those for the first scenario in the slab sampling section on the IEEE 123-bus
feeder, and the runtime is 28 seconds.
}
\rev{
\section{Comparison between Distributed Tensor Completion and Decentralized Matrix Completion}
In \cite{Sagan2019}, a decentralized model-based matrix completion method was developed by building  the data matrix using multiple time-step
data.
We compare our proposed distributed tensor completion approach to the decentralized model-based matrix completion method for Case 1  in Sec.~\ref{sec:slabsampling} and the comparison results are shown in Table~\ref{tab:comparisonMC}.}

\rev{
\begin{table}
	\renewcommand{\arraystretch}{1}
	\captionsetup{labelfont={color=black},font={color=black}}
	\caption{Comparison of the proposed algorithm and decentralized matrix completion }\vspace{-10pt}
\label{tab:comparisonMC}
	\begin{center}
		\begin{tabular}{l c c c c } 
			\toprule
	            & \res{MAPE($|V|$)}  & \res{MAE($\theta$)} & \res{MAE($P$)} & \res{MAE($Q$)} \\
			\midrule
			\res{TC} &  \res{$0.5922\%$} & \res{$0.7758$} & \res{$0.0279$} & \res{$0.0105$} \\
		\res{MC} & \res{$0.4263\%$} & \res{$0.2134$} &\res{$0.0136$} & \res{$0.0068$} \\
			\bottomrule
		\end{tabular}\vspace{-15pt}
	\end{center}
\end{table}
}
\rev{The results for
the decentralized matrix completion are generally better than our model but
of the same order for MAPE for the voltage magnitude, and MAE of the voltage angle and active power. Note that our proposed approach achieves
this result without the use of the model information required
by the decentralized matrix completion method \cite{Sagan2019}. Also the runtime for our proposed method is 16 seconds, and that for the matrix completion is around 10 minutes. But it is worth mentioning that the decentralized matrix completion model can perform well enough at the same percentage of measurements  with only one time step data  due to the power flow constraints.
}

\bibliographystyle{IEEEtran}
\bibliography{IEEEabrv,Bibliography,references}

\vfill


\end{document}